\documentclass[12pt]{article}
\usepackage{cite,amsfonts,amsmath}
\usepackage{epic}

\renewcommand{\appendix}{%
\renewcommand{\section}{%
\newpage\thispagestyle{plain}%
\secdef\Appendix\sAppendix}%
\setcounter{section}{0}%
\renewcommand{\thesection}{\Alph{section}}%
}
\newcommand{\Appendix}[2][?]{%
\refstepcounter{section}%
\addcontentsline{toc}{appendix}%
{\protect\numberline{\appendixname~\thesection}#1}%
{\flushleft\LARGE\bfseries\appendixname\ \thesection\par
\centering#2\par}%
\sectionmark{#1}\vspace{\baselineskip}}

\newcommand{\sAppendix}[1]{%
{\flushright\large\bfseries\appendixname\par
\centering#1\par}%
\vspace{\baselineskip}}

\def\be{\begin{equation}}

\def\bea{\begin{eqnarray}}

\def\eea{\end{eqnarray}}

\begin{document}

\pagestyle{empty}

\begin{center}

\textsf{\Huge {\bf A nested sequence of projectors (2): Multiparameter multistate statistical models,
Hamiltonians, $S$-matrices}}

\vspace{10mm}

{\large B. Abdesselam$^{a,}$\footnote{This work was carried out at the Centre de Physique Th\'eorique de l'Ecole
Polytechnique, Palaiseau.}$^{,}$\footnote{Email: boucif@cpht.polytechnique.fr and  boucif@yahoo.fr} and
A. Chakrabarti$^{b,}$\footnote{Email: chakra@cpht.polytechnique.fr}}

 \vspace{5mm}

 \emph{$^a$ Laboratoire de Physique Quantique de la
Mati\`ere et de Mod\'elisations Math\'ematiques, Centre
Universitaire de Mascara, 29000-Mascara, Alg\'erie\\
and \\
Laboratoire de Physique Th\'eorique, Universit\'e d'Oran
Es-S\'enia, 31100-Oran, Alg\'erie}
 \\
 \vspace{3mm}
 \emph{$^b$ Centre de Physique Th{\'e}orique, CNRS UMR 7644}
 \\
 \emph{Ecole Polytechnique, 91128 Palaiseau Cedex, France.}
 \vspace{3mm}

\end{center}

\begin{abstract}
Our starting point is a class of braid matrices, presented in a
previous paper, constructed on a basis of a nested sequence of
projectors. Statistical models associated to such $N^2\times N^2$
matrices for odd $N$ are studied here. Presence of $\frac
12\left(N+3\right)\left(N-1\right)$ free parameters is the crucial
feature of our models, setting them apart from other well-known
ones. There are $N$ possible states at each site. The trace of the transfer matrix is shown to depend
on $\frac 12\left(N-1\right)$ parameters. For order $r$, $N$ eigenvalues consitute the trace and the remaining
$\left(N^r-N\right)$ eigenvalues involving the full range
of parameters come in zero-sum multiplets formed
by the $r$-th roots of unity, or lower dimensional multiplets
corresponding to factors of the order $r$ when $r$ is not a prime
number. The modulus of any eigenvalue is of the form
$e^{\mu\theta}$, where $\mu$ is a linear combination of the free
parameters, $\theta$ being the spectral parameter. For $r$ a prime
number an amusing relation of the number of multiplets with a
theorem of Fermat is pointed out. Chain Hamiltonians and
potentials corresponding to factorizable $S$-matrices are
constructed starting from our braid matrices. Perspectives are
discussed.
\end{abstract}

\rightline{math.QA/0601584}
\rightline{January 2006}

\pagestyle{plain}
\setcounter{page}{1}

%\tableofcontents

\newpage
\section{Introduction}
\setcounter{equation}{0}

The most salient feature of the class of braid matrices presented in ref. \cite{1}, setting it apart from other known
examples, is the number of free parameters. This class was obtained for $N^2\times N^2$ braid matrices for {\it odd}
\begin{equation}
N=\left(2p-1\right),\qquad \left(p=1,2,\ldots\right).
\end{equation}
Such matrices, depending on a spectral parameter $\theta$ and satisfying the braid equation
\begin{equation}
\widehat{R}_{12}\left(\theta-\theta'\right)\widehat{R}_{23}\left(\theta\right)\widehat{R}_{12}\left(\theta'\right)
=\widehat{R}_{23}\left(\theta'\right)\widehat{R}_{12}\left(\theta\right)\widehat{R}_{23}\left(\theta-\theta'\right)
\end{equation}
have $\frac 12\left(N+3\right)\left(N-1\right)$ free parameters when the overall normalization is fixed. Thus for
$N=3,\,5,\,7,\ldots$, the respective number of parameters are $6,\,16,\,30,\ldots$. These parameters appear in the
coefficients of the $N^2$ projectors (the "{\em nested sequence}" defined in ref. \cite{1}) providing the basis of
$\widehat{R}\left(\theta\right)$.

The projectors are defined as follows. Let $\left(ij\right)$ be the $N\times N$ matrix with a single non-zero element
1 on row $i$ and column $j$. Then $N^2$ projectors are defined, with $\epsilon=\pm$ and $\overline{i}=N-i+1$, as
\begin{eqnarray}
&& P_{pp}=(pp)\otimes (pp),\nonumber\\
&&2P_{pi\left(\epsilon\right)}=(pp)\otimes\left[(ii)+(\overline{i}\overline{i})+\epsilon\left((i\overline{i})+
(\overline{i}i)\right)\right],\nonumber\\
&&2P_{ip\left(\epsilon\right)}=\left[(ii)+(\overline{i}\overline{i})+\epsilon\left((i\overline{i})+
(\overline{i}i)\right)\right]\otimes(pp),\nonumber\\
&&2P_{ij\left(\epsilon\right)}=(ii)\otimes(jj)+(\overline{i}\overline{i})\otimes(\overline{j}\overline{j})+
\epsilon\left[(i\overline{i})\otimes(j\overline{j})+(\overline{i}i)\otimes(\overline{j}j)\right],\nonumber\\
&&2P_{i\overline{j}\left(\epsilon\right)}=(ii)\otimes(\overline{j}\overline{j})+(\overline{i}\overline{i})\otimes
(jj)+\epsilon\left[(i\overline{i})\otimes(\overline{j}j)+(\overline{i}i)\otimes(j\overline{j})\right],
\end{eqnarray}
where, from (1.1),
\begin{equation}
i=1,\,2,\ldots,\,p-1,\qquad
\overline{i}=N-i+1=2p-1,\,2p-2,\ldots,\,p+1,\qquad
p=\frac 12\left(N+1\right).
\end{equation}
The projectors satisfy (with $\alpha$, $\beta$ standing for triplets $\left(i,j,\epsilon\right)$)
\begin{equation}
P_\alpha P_\beta=\delta_{\alpha\beta}P_{\alpha},\qquad \sum_\alpha P_\alpha^2=I_{N^2\times N^2}.
\end{equation}
Their total number is
\begin{equation}
1+4\left(p-1\right)+4\left(p-1\right)^2=\left(2p-1\right)^2=N^2.
\end{equation}
For our class of solutions, normalizing to 1 the coefficient of $P_{pp}$,
\begin{equation}
\widehat{R}\left(\theta\right)=P_{pp}+\sum_{i,\epsilon}\left(e^{m_{pi}^{(\epsilon)}\theta}P_{pi(\epsilon)}
+e^{m_{ip}^{(\epsilon)}\theta}P_{ip(\epsilon)}\right)+\sum_{i,j,\epsilon}\left(e^{m_{ij}^{(\epsilon)}\theta}P_{ij(\epsilon)}
+e^{m_{i\overline{j}}^{(\epsilon)}\theta}P_{i\overline{j}(\epsilon)}\right),
\end{equation}
with the crucial constraint
\begin{equation}
m_{ij}^{(\epsilon)}=m_{i\overline{j}}^{(\epsilon)},\qquad \left(\overline{j}=N-j+1=2p-j\right).
\end{equation}
This sufficient and necessary constraint concerning the coefficient of $\theta$ in the exponents, leaves
\begin{equation}
\frac 12\left(N+3\right)\left(N-1\right)
\end{equation}
free parameters. For $N=3$ one thus obtains, with 6 free parameters,
\begin{equation}
\widehat{R}\left(\theta\right)=\left|\begin{array}{ccccccccc}
  a_+ & 0 & 0 & 0 & 0 & 0 & 0 & 0 & a_- \\
  0 & b_+ & 0 & 0 & 0 & 0 & 0 & b_- & 0 \\
  0 & 0 & a_+ & 0 & 0 & 0 & a_- & 0 & 0 \\
  0 & 0 & 0 & c_+ & 0 & c_- & 0 & 0 & 0 \\
  0 & 0 & 0 & 0 & 1 & 0 & 0 & 0 & 0 \\
  0 & 0 & 0 & c_- & 0 & c_+ & 0 & 0 & 0 \\
  0 & 0 & a_- & 0 & 0 & 0 & a_+ & 0 & 0 \\
  0 & b_- & 0 & 0 & 0 & 0 & 0 & b_+ & 0 \\
  a_- & 0 & 0 & 0 & 0 & 0 & 0 & 0 & a_+ \\
\end{array}\right|
\end{equation}
where
\begin{equation}
a_{\pm}=\frac 12\left(e^{m_{11}^{(+)}\theta}\pm e^{m_{11}^{(-)}\theta}\right),\qquad b_{\pm}=
\frac 12\left(e^{m_{12}^{(+)}\theta}\pm e^{m_{12}^{(-)}\theta}\right),\qquad c_{\pm}=
\frac 12\left(e^{m_{21}^{(+)}\theta}\pm e^{m_{21}^{(-)}\theta}\right).
\end{equation}
The parameters $a_{\pm}$ are each repeated in (1.10) according to (1.8), since
\begin{equation}
m_{1\overline{1}}^{(\pm)}=m_{11}^{(\pm)},\qquad \left(\overline{1}=3\right).
\end{equation}
This is the case we will study mostly in the following sections. The corresponding results for $N>3$ will be indicated
briefly. For example, the generalization of the considerations below in this section for $N>3$ is entirely
straight-forward. To explore the statistical model associated to (1.10) one starts by constructing explicit
representations of the monodromy matrices $t_{ij}^{(r)}\left(\theta\right)$ of successive orders $\left(r=1,2,3,\ldots
\right)$ obtained by taking coproducts of the fundamental $3\times 3$ blocks (with the same $\theta$ for each
factor)
\begin{equation}
t_{ij}^{(r)}=\sum_{j_1,\ldots,j_{r-1}}t_{ij_1}\otimes t_{j_1j_2}\otimes\cdots\otimes t_{j_{r-1},j}.
\end{equation}
For $N=3$,
\begin{equation}
t^{(r)}=\left|\begin{array}{ccccc}
t_{11}^{(r)}& & t_{12}^{(r)} & & t_{1\overline{1}}^{(r)} \\
&&&&\\
t_{21}^{(r)} && t_{22}^{(r)} && t_{2\overline{1}}^{(r)} \\
&&&&\\
t_{\overline{1}1}^{(r)} && t_{\overline{1}2}^{(r)} && t_{\overline{1}\overline{1}}^{(r)} \\
\end{array}\right|.
\end{equation}
If the $\widehat{R}tt$ equation for the blocks $t_{ij}^{(r)}$ (appendix C),
\begin{equation}
\widehat{R}\left(\theta-\theta'\right)\left(t^{(r)}\left(\theta\right)\otimes
t^{(r)}\left(\theta'\right)\right)=\left(t^{(r)}\left(\theta'\right)\otimes
t^{(r)}\left(\theta\right)\right)\widehat{R}\left(\theta-\theta'\right)
\end{equation}
is satisfied for $r=1$, then the coproduct construction (1.13) ensures that (1.15) is satisfied for all higher values
$r=2,\,3,\ldots$. The solution for $r=1$ is given by
\begin{equation}
t^{(1)}\left(\theta\right)\equiv t\left(\theta\right)=P\widehat{R}\left(\theta\right)=R\left(\theta\right),
\end{equation}
where $P$ is the permutation matrix
\begin{equation}
P=\sum_{ij}\left(ij\right)\otimes\left(ji\right)
\end{equation}
and $R\left(\theta\right)$ is the Yang-Baxter (YB) matrix. This is a standard result valid generally for solutions
of (1.2). (See appendix B of ref. \cite{1} for sources cited.)

The transfer matrix, for each order $r$, is defined to be the trace (with argument $\theta$)
\begin{equation}
T^{(r)}=t_{11}^{(r)}+t_{22}^{(r)}+t_{\overline{1}\overline{1}}^{(r)}.
\end{equation}
The properties of the model depend crucially on the eigenvalues of $T^{(r)}$. Refs. \cite{2,3,4} provide ample
information citing numerous basic sources.

So our basic task will be to construct the eigenstates and eigenvalues of $T^{(r)}\left(\theta\right)$. Remarkable
feature following from (1.10) (and more generally from (1.7)) will be presented in the following sections and appendices.

We will also construct chain Hamiltonians and potentials leading to factorizable $S$-matrices starting from our class of
$\widehat{R}\left(\theta\right)$.

Concerning each aspect we will try to display the role of our multiple parameters. For all
\begin{equation}
m_{ij}^{(+)}\theta\geq m_{ij}^{(-)}\theta
\end{equation}
the elements of $\widehat{R}\left(\theta\right)$and hence the Boltzmann weights are non-negative, consistent with
physical interpretations. For {\it definiteness} we consider the sector, say
\begin{equation}
m_{11}^{(+)}> m_{11}^{(-)}>m_{12}^{(+)}> m_{12}^{(-)}m_{21}^{(+)}>
m_{21}^{(-)},\qquad \theta\geq 0
\end{equation}
of (1.10). The eigenvalues will be ordered differently for other
sectors. They can be considered separately.

\section{Transfer matrix, eigenvectors, eigenvalues ($N=3$): crucial features}
\setcounter{equation}{0}

We start by signalling some crucial features to be encountered below in the explicit constructions restricted (in this
section) to $N=3$.
\begin{enumerate}
\item The trace of the transfer matrix (1.17) of order $r$ will  turn out to be
\begin{equation}
\hbox{tr}\left(T^{(r)}\left(\theta\right)\right)=2e^{rm_{11}^{(+)}\theta}+1.
\end{equation}
Of the six parameters $\left(m_{11}^{(\pm)},m_{12}^{(\pm)},m_{21}^{(\pm)}\right)$ of (1.11) {\it only} $m_{11}^{(+)}$
appears in the trace. A simple explanation of this fact will be given after discussing the generalization for $N>3$.

\item The eigenvalue $e^{rm_{11}^{(+)}\theta}$ is obtained exactly {\it twice} for each $r$ and the value 1 only
{\it once}.

\item The remaining $\left(3^r-3\right)$ eigenvalues occur in multiplets of {\it zero sum} due to the presence of
{\it roots of unity}. Hence they do not contribute to the trace. For $r$ a prime number there will be "$r$-plets"
(and possibly "$nr$-plets", $n$ being an integer)
\begin{equation}
e^{\mu\theta}\left(1,e^{\frac{2\pi i}{r}},e^{\frac{2\pi i}{r}\cdot 2},\ldots,,e^{\frac{2\pi i}{r}\cdot\left(r-1\right)}
\right),
\end{equation}
where $\mu$ is a {\it linear} combination of the parameters $m_{ij}^{(\pm)}$. When $r$ is factorizable lower order
multiplets can be present corresponding to the factors. Thus for $r=4$ one obtains both doublets and quadruplets
\begin{equation}
e^{\mu_2\theta}\left(1,-1\right),\qquad e^{\mu_4\theta}\left(1,i,-1,-i\right),
\end{equation}
with appropriate linear combinations $\mu_2$, $\mu_4$ to be displayed below.

\item Apart from possible roots of unity phase factors the modulus of each eigenvalue is a simple exponential of the
type $e^{\mu\theta}$ of (2.2). For $r=3$, for example, one obtains for $\mu$ the values (appendix A)
\begin{eqnarray}
&&3m_{11}^{(+)},\,\left(m_{11}^{(+)}+2m_{11}^{(-)}\right),\nonumber\\
&&\left(m_{11}^{(+)}+m_{12}^{(+)}+m_{21}^{(+)}\right),\,\left(m_{11}^{(+)}+m_{12}^{(-)}+m_{21}^{(-)}\right),\nonumber\\
&&\left(m_{11}^{(-)}+m_{12}^{(+)}+m_{21}^{(-)}\right),\,\left(m_{11}^{(-)}+m_{12}^{(-)}+m_{21}^{(+)}\right),\nonumber\\
&&\left(m_{12}^{(+)}+m_{21}^{(+)}\right),\,\left(m_{12}^{(-)}+m_{21}^{(-)}\right),\nonumber\\
&&0.
\end{eqnarray}
Along with roots of unity factors these provides all the 27 eigenvalues as will be shown below.

\item The values of $\mu$ depend crucially on the subspaces, to be introduced below, which are invariant under the
action of $T^{(r)}\left(\theta\right)$, the transfer matrix.
\end{enumerate}

\subsection{Construction of $T^{(r)}\left(\theta\right)$ for $N=3$}

The standard construction of the fundamental $3\times 3$ block matrices $t_{ij}\left(\theta\right)$ implementing (1.15),
(1.16), (1.10), (1.11) leads to (for $t^{(1)}\left(\theta\right)\equiv t\left(\theta\right)$ with $\overline{1}=3$)
\begin{eqnarray}
&&t_{11}\left(\theta\right)=\left|\begin{array}{ccc}
  a_+ & 0 & 0 \\
  0 & 0 & 0 \\
  0 & 0 & a_- \\
\end{array}\right|,\qquad t_{12}\left(\theta\right)=\left|\begin{array}{ccc}
  0 & 0 & 0 \\
  c_+ & 0 & c_- \\
  0 & 0 & 0 \\
\end{array}\right|,\qquad t_{1\overline{1}}\left(\theta\right)=\left|\begin{array}{ccc}
  0 & 0 & a_- \\
  0 & 0 & 0 \\
  a_+ & 0 & 0 \\
\end{array}\right|,\nonumber\\
&&t_{21}\left(\theta\right)=\left|\begin{array}{ccc}
  0 & b_+ & 0 \\
  0 & 0 & 0 \\
  0 & b_- & 0 \\
\end{array}\right|,\qquad t_{22}\left(\theta\right)=\left|\begin{array}{ccc}
  0 & 0 & 0 \\
  0 & 1 & 0 \\
  0 & 0 & 0 \\
\end{array}\right|,\qquad t_{2\overline{1}}\left(\theta\right)=\left|\begin{array}{ccc}
  0 & b_- & 0 \\
  0 & 0 & 0 \\
  0 & b_+ & 0 \\
\end{array}\right|,\\
&&t_{\overline{1}1}\left(\theta\right)=\left|\begin{array}{ccc}
  0 & 0 & a_+ \\
  0 & 0 & 0 \\
  a_- & 0 & 0 \\
\end{array}\right|,\qquad t_{\overline{1}2}\left(\theta\right)=\left|\begin{array}{ccc}
  0 & 0 & 0 \\
  c_- & 0 & c_+ \\
  0 & 0 & 0 \\
\end{array}\right|,\qquad t_{\overline{1}\overline{1}}\left(\theta\right)=\left|\begin{array}{ccc}
  a_- & 0 & 0 \\
  0 & 0 & 0 \\
  0 & 0 & a_+ \\
\end{array}\right|,\nonumber
\end{eqnarray}
where, from (1.11),
\begin{equation}
\left(a_+\pm a_-\right)=e^{m_{11}^{(\pm)}\theta},\qquad \left(b_+\pm b_-\right)=e^{m_{12}^{(\pm)}\theta},\qquad
\left(c_+\pm c_-\right)=e^{m_{21}^{(\pm)}\theta}.
\end{equation}
One now has to be implement these in (1.13), (1.14) and (1.18) to obtain $T^{(r)}\left(\theta\right)$. Then one proceeds
to construct eigenvalues of $T^{(r)}\left(\theta\right)$.

\subsection{Subspaces invariant under the action of $T^{(r)}\left(\theta\right)$}

We start by introducing convenient, compact notations. The state vectors of the fundamental representation (2.5) are
denoted as
\begin{equation}
\left(\left|\begin{array}{c} 1 \\ 0 \\  0 \\\end{array}\right\rangle,\,\left|\begin{array}{c}
  0 \\  1 \\  0 \\\end{array}\right\rangle,\,\left|\begin{array}{c}  0 \\  0 \\ 1 \\
\end{array}\right\rangle\right)\equiv\left(\left|\begin{array}{c}
  1 \\\end{array}\right\rangle,\,\left|\begin{array}{c}
  2 \\\end{array}\right\rangle,\,\left|\begin{array}{c}
  \overline{1} \\\end{array}\right\rangle\right).
\end{equation}
Tensor products for higher orders are constructed as
\begin{equation}
\left(\left|1\right\rangle\otimes\left|1\right\rangle,\left|1\right\rangle\otimes\left|2\right\rangle,
\left|1\right\rangle\otimes\left|\overline{1}\right\rangle,\ldots\right)\equiv\left(\left|11\right\rangle,
\left|12\right\rangle,\left|1\overline{1}\right\rangle,\ldots\right)
\end{equation}
and so on in evident continuation. The {\it order} of the labels $(1,2,\overline{1})$ will indicate the tensor product
structure. Thus, for example,
\begin{equation}
\left|1\right\rangle\otimes\left|1\right\rangle\otimes\left|2\right\rangle\otimes
\left|\overline{1}\right\rangle\otimes\left|1\right\rangle\equiv\left|112\overline{1}2\right\rangle
\end{equation}

Corresponding to the $r$-th order coproduct, $T^{(r)}\left(\theta\right)$ acts on a space
spanned by $3^r$ states (for $N=3$). Let
\begin{equation}
S\left(r,k\right),\qquad \left(k=0,1,\ldots,r\right)
\end{equation}
denote the subspaces labeled by $k$, the {\it multiplicity} of the index $2$. The coefficients of different
power of $x$ in the expansion
\begin{equation}
\left(x+2\right)^r=1\cdot x^r+2rx^{r-1}+\cdots+2^{r-k}\binom r{r-k}x^k+\cdots+2^r
\end{equation}
give the number of states in the respective subspaces. Setting $x=1$ one obtains the total number
\begin{equation}
\left(1+2\right)^r=3^r.
\end{equation}
For example, for $r=3$, one obtains the subspaces,
\begin{eqnarray}
&&S\left(3,3\right):\qquad \left|222\right\rangle\nonumber\\
&&S\left(3,2\right):\qquad \left|221\right\rangle,\left|22\overline{1}\right\rangle,\left|212\right\rangle,
\left|2\overline{1}2\right\rangle,\left|122\right\rangle,\left|\overline{1}22\right\rangle,\nonumber\\
&&S\left(3,1\right):\qquad \left|211\right\rangle,\left|21\overline{1}\right\rangle,\left|2\overline{1}1\right\rangle,
\left|2\overline{1}\overline{1}\right\rangle\nonumber\\
&&\phantom{S\left(3,1\right):\qquad}\left|121\right\rangle,\left|12\overline{1}\right\rangle,
\left|\overline{1}21\right\rangle,
\left|\overline{1}2\overline{1}\right\rangle\nonumber\\
&&\phantom{S\left(3,1\right):\qquad}\left|112\right\rangle,\left|1\overline{1}2\right\rangle,
\left|\overline{1}12\right\rangle,\left|\overline{1}\overline{1}2\right\rangle\nonumber\\
&&S\left(3,0\right):\qquad \left|111\right\rangle,\left|11\overline{1}\right\rangle,\left|1\overline{1}1\right\rangle,
\left|\overline{1}11\right\rangle\nonumber\\
&&\phantom{S\left(3,0\right):\qquad}\left|\overline{1}\overline{1}\overline{1}\right\rangle,
\left|\overline{1}\overline{1}1\right\rangle,\left|\overline{1}1\overline{1}\right\rangle,
\left|1\overline{1}\overline{1}\right\rangle
\end{eqnarray}
A striking and most helpful consequence of the structure of the matrices (2.5) and their coproducts is: {\it each subspace
$S\left(r,k\right)$ is invariant under the action of $T^{(r)}\left(\theta\right)$}. This facilitates considerably the
construction of eigenstates. One works on lower dimensional spaces.

One possible approach is as follows: One selects any one state from the $2^{r-k}\binom r{r-k}$ states of
$S\left(r,k\right)$ and computes the action of $T^{(r)}\left(\theta\right)$ on it. One gets on the r.h.s. a linear
combination of states belonging to $S\left(r,k\right)$. Thus, for example,
\begin{eqnarray}
&&T^{(4)}\left(\theta\right)\left|1111\right\rangle=\left(a_+^4+a_-^4\right)\left|1111\right\rangle+
2a_+^2a_-^2\left(\left|\overline{1}\overline{1}\overline{1}\overline{1}\right\rangle+
\left|1\overline{1}1\overline{1}\right\rangle+\left|\overline{1}1\overline{1}1\right\rangle\right)+\nonumber\\
&&\phantom{T^{(4)}\left(\theta\right)\left|1111\right\rangle=}\left(a_+^2+a_-^2\right)a_+a_-\left(\left|11\overline{1}\overline{1}\right\rangle+
\left|\overline{1}\overline{1}11\right\rangle+\left|\overline{1}11\overline{1}\right\rangle+
\left|1\overline{1}\overline{1}1\right\rangle\right),
\end{eqnarray}
where $a_{\pm}=\frac 12\left(e^{m_{11}^{(+)}\theta}\pm e^{m_{11}^{(-)}\theta}\right)$ as noted before. Next one
computes the action of $T^{(4)}$ successively on the other states appearing on the right. This continues until
one obtains the coefficients for a closed subsystem. Then one searches for linear combinations such that under
$T^{(4)}$ it is reproduced to within a factor. Thus one systematically obtains all eigenvectors and eigenvalues, for
the subspace $S\left(r,k\right)$. For our class one has to solve systems of {\it linear} equations with fairly simple
coefficient. Even the 81 eigenstates and eigenvalues for $r=4$ were obtained directly without using a computer program
and without any real difficulties.

We have thus obtained exhaustive solutions for $r=1,\,2,\,3,\,4$. The corresponding $3,\,9,\,27$ and $81$ eigenvalues
are presented in appendix A. We have also obtained explicitly all the corresponding eigenstates. For brevity they are
not presented here. The eigenvalues of the appendix A fully illustrate the crucial properties (1) to (5) signalled at
the start of this section. In the following section we indicate a related but somewhat differently formulated approach
for various comparisons.

\section{Linear constraints for eigenvectors for $N=3$ and comparison with algebraic Bethe ansatz}
\setcounter{equation}{0}

In section 2 we have noted how, exploiting the invariance of the subspaces $S\left(r,k\right)$ defined by (2.10)
one can construct step by step all the eigenstates. The comments following (2.14) indicate how the relevant linear
equations are obtained. We formulate below the approach in a systematic, explicit fashion.

Starting with (2.5) and (2.6) for $t_{ij}^{(1)}\left(\theta\right)=t_{ij}\left(\theta\right)$ we define the operators
\begin{eqnarray}
&&U=b_+\left(-\theta\right)t_{21}\left(\theta\right)+b_-\left(-\theta\right)t_{2\overline{1}}\left(\theta\right)
=\left(\begin{array}{ccc}
  0 & 1 & 0 \\
  0 & 0 & 0 \\
  0 & 0 & 0 \\
\end{array}\right),\nonumber\\
&&A=t_{22}\left(\theta\right)=\left(\begin{array}{ccc}
  0 & 0 & 0 \\
  0 & 1 & 0 \\
  0 & 0 & 0 \\
\end{array}\right),\nonumber\\
&&D=b_-\left(-\theta\right)t_{21}\left(\theta\right)+b_+\left(-\theta\right)t_{2\overline{1}}\left(\theta\right)
=\left(\begin{array}{ccc}
  0 & 0 & 0 \\
  0 & 0 & 0 \\
  0 & 1 & 0 \\
\end{array}\right).
\end{eqnarray}
Then (suppressing arguments $\theta$ of $t_{ij}$)
\begin{eqnarray}
&&t_{11}\left(A,U,D\right)=\left(0,a_+U,a_-D\right)\nonumber\\
&&t_{12}\left(A,U,D\right)=\left(0,c_+A,c_-A\right)\nonumber\\
&&t_{1\overline{1}}\left(A,U,D\right)=\left(0,a_+D,a_-U\right)\nonumber\\
&&t_{21}\left(A,U,D\right)=\left(b_+U+b_-D,0,0\right)\nonumber\\
&&t_{22}\left(A,U,D\right)=\left(A,0,0\right)\nonumber\\
&&t_{2\overline{1}}\left(A,U,D\right)=\left(b_-U+b_+D,0,0\right)\nonumber\\
&&t_{\overline{1}1}\left(A,U,D\right)=\left(0,a_-D,a_+U\right)\nonumber\\
&&t_{\overline{1}2}\left(A,U,D\right)=\left(0,c_-A,c_+A\right)\nonumber\\
&&t_{\overline{1}\overline{1}}\left(A,U,D\right)=\left(0,a_-U,a_+D\right).
\end{eqnarray}
Also from (2.7) and (3.1)
\begin{equation}
U\left|2\right\rangle=\left|1\right\rangle,\qquad A\left|2\right\rangle=\left|2\right\rangle,\qquad
D\left|2\right\rangle=\left|\overline{1}\right\rangle.
\end{equation}
For any $r$, starting with $S\left(r,r\right)$ one obtains the basic eigenstate (trivially since  $S\left(r,r\right)$
is of dimension 1),
\begin{equation}
T^{(r)}\left(\theta\right)\left|22\ldots 2\right\rangle=1\left|22\ldots 2\right\rangle.
\end{equation}
Now one moves up in $\left(r-k\right)$ stepwise.

\paragraph{$\underline{S\left(r,r-1\right)\,(\dim 2r}$):} With $2r$ coefficients $\left(u_i,d_i\right)$
$\left(i=1,\ldots,r\right)$ one can label the states as
\begin{eqnarray}
\left(\begin{array}{l}
  \left(u_1U+d_1D\right)\otimes A\otimes\cdots\otimes A \\
  +A\otimes\left(u_2U+d_2D\right)\otimes A\otimes\cdots\otimes A \\
  \vdots \\
  + A\otimes A\otimes\cdots\otimes A\otimes\left(u_rU+d_rD\right)\\
\end{array}\right)\left|22\ldots 2\right\rangle.
\end{eqnarray}
The action of $T^{(r)}\left(\theta\right)$ on these leads to a linear system of equations in $\left(u_i,d_i\right)$
corresponding to eigenstates. For $S\left(r,r-1\right)$ the solution is particularly simple. Define
\begin{equation}
\left|\omega,\epsilon\right\rangle=\left(\begin{array}{l}
  A\otimes A\otimes\cdots\otimes A\otimes\left(U+\epsilon D\right) \\
  +\omega A\otimes A\otimes\cdots\otimes\left(U+\epsilon D\right)\otimes A \\
  +\omega^2 A\otimes A\otimes\cdots\otimes\left(U+\epsilon D\right)\otimes A\otimes A \\
  \vdots \\
  +\omega^{r-1}\left(U+\epsilon D\right)\otimes A\otimes A\otimes\cdots\otimes \\
\end{array}\right)
\left|22\ldots 2\right\rangle,
\end{equation}
where $\epsilon=\pm$ and $\omega$ can have $r$ values (as a $r$-th root of unity)
\begin{equation}
\omega=\left(1,e^{\frac{i2\pi}{r}},\ldots,e^{\frac{i2\pi}{r}\cdot\left(r-1\right)}\right).
\end{equation}
One obtains
\begin{equation}
T^{(r)}\left(\theta\right)\left|\omega,\epsilon\right\rangle=\omega^{r-1}e^{(m_{12}^{(\epsilon)}+
m_{21}^{(\epsilon)})\theta}\left|\omega,\epsilon\right\rangle.
\end{equation}
The $2r$ eigenvalues are
\begin{equation}
e^{(m_{12}^{(\epsilon)}+
m_{21}^{(\epsilon)})\theta}\left(1,e^{\frac{i2\pi}{r}},\ldots,e^{\frac{i2\pi}{r}\left(r-1\right)}\right).
\end{equation}
The next step is $k=r-2$.

\paragraph{$\underline{S\left(r,r-2\right)\,(\dim 2r\left(r-1\right))}$:} A set of states spanning this subspace
is given by
\begin{equation}
\sum_{i\neq j}\left(A\otimes\cdots\otimes A\otimes\left(u_iU+d_i D\right)\otimes A\otimes\cdots\otimes
\left(u_jU+d_j D\right)\otimes A\otimes\cdots\otimes A\right)
\left|22\ldots 2\right\rangle.
\end{equation}
The parameters $\left(u,d\right)$ have to be constrained to obtain eigenstates. At each step one obtains sets of
linear constraints. The pattern is now evident. At each step one inserts, as in (3.10), $\left(r-k\right)$ factors
of the type $\left(u_iU+d_i D\right)$ excluding their coincidence.

Finally for $k=0$, one has $S\left(r,0\right)$ of dimension $2^r$. Here a basis spanning the subspace can be labeled
as
\begin{equation}
\left(\sum_{i}\left(u^{(i)}_1U+d^{(i)}_1 D\right)\otimes \left(u^{(i)}_2U+d^{(i)}_2 D\right)\otimes\cdots\otimes
\left(u^{(i)}_rU+d^{(i)}_r D\right)\right)
\left|22\ldots 2\right\rangle.
\end{equation}
Since there are no label $\left|2\right\rangle$ left, it can be shown from (3.2) that $T^{(r)}\left(\theta\right)$,
acting on $S\left(r,0\right)$ simplifies to
\begin{eqnarray}
&&T^{(r)}\approx t_{11}^{(r)}+t_{\overline{1}\overline{1}}^{(r)}\approx
t_{11}\otimes t_{11}^{(r-1)}+t_{1\overline{1}}\otimes t_{\overline{1}1}^{(r-1)}+
t_{\overline{1}1}\otimes t_{1\overline{1}}^{(r-1)}+t_{\overline{1}\overline{1}}
\otimes t_{\overline{1}\overline{1}}^{(r-1)}.
\end{eqnarray}
In successive steps ($t^{(r-1)}\rightarrow t\otimes t^{(r-2)}$ and so on) only the indices $\left(1,\overline{1}\right)$
need be retained. They only give non zero contributions. From (3.2)
\begin{eqnarray}
&&t_{11}\left(u_s^{(i)}U+d_s^{(i)}D\right)=\left(a_+u_s^{(i)}U+a_-d_s^{(i)}D\right)\equiv X_{11}^{(i)}\left(s\right),
\nonumber\\
&&t_{1\overline{1}}\left(u_s^{(i)}U+d_s^{(i)}D\right)=\left(a_-d_s^{(i)}U+a_+u_s^{(i)}D\right)\equiv
X_{1\overline{1}}^{(i)}\left(s\right),\nonumber\\
&&t_{\overline{1}1}\left(u_s^{(i)}U+d_s^{(i)}D\right)=\left(a_+d_s^{(i)}U+a_-u_s^{(i)}D\right)\equiv
X_{\overline{1}1}^{(i)}\left(s\right),\nonumber\\
&&t_{\overline{1}\overline{1}}\left(u_s^{(i)}U+d_s^{(i)}D\right)=\left(a_-u_s^{(i)}U+a_+d_s^{(i)}D\right)\equiv
X_{\overline{1}\overline{1}}^{(i)}\left(s\right).
\end{eqnarray}

Define, with indices taking values $\left(1,\overline{1}\right)$ only,
\begin{equation}
X_{ab}^{(i)}\left(1,2,\ldots,r\right)=\sum_{b_1,\ldots,b_{r-1}}X_{ab_1}^{(i)}\left(1\right)\otimes
X_{b_1b_2}^{(i)}\left(2\right)\otimes\cdots\otimes X_{b_{r-1}b}^{(i)}\left(r\right).
\end{equation}
The action of $T^{(r)}\left(\theta\right)$ on the generic state (3.11) finally reduces to
\begin{equation}
\left(\sum_{i=1}^r\left(X_{11}^{(i)}\left(1,2,\ldots,r\right)+X_{\overline{1}\overline{1}}^{(i)}\left(1,2,\ldots,r\right)
\right)\right)\left|22\ldots 2\right\rangle.
\end{equation}
It is of particular interest to see what parameterizations in (3.11) corresponds to the {\it two} eigenstates (and two
only for any $r$) that contribute to the trace.

For $r=2$,
\begin{eqnarray}
&&T^{(2)}\left(\theta\right)\left(\left|11\right\rangle+\left|\overline{1}\overline{1}\right\rangle\right)=e^{2m_{11}^{(+)}\theta}
\left(\left|11\right\rangle+\left|\overline{1}\overline{1}\right\rangle\right),\nonumber\\
&&T^{(2)}\left(\theta\right)\left(\left|1\overline{1}\right\rangle+\left|\overline{1}1\right\rangle\right)=e^{2m_{11}^{(+)}\theta}
\left(\left|1\overline{1}\right\rangle+\left|\overline{1}1\right\rangle\right).
\end{eqnarray}
Along with
\begin{eqnarray}
&&T^{(2)}\left(\theta\right)\left|22\right\rangle=\left|22\right\rangle
\end{eqnarray}
the two states of (3.16) yield
\begin{eqnarray}
&&\hbox{tr}\left(T^{(2)}\left(\theta\right)\right)=2e^{2m_{11}^{(+)}\theta}+1.
\end{eqnarray}
For other six states provide the $3$ zero sum doublets of (A.4), (A.5), (A.6).

For $r=3$. Apart from $8$ zero sum triplets of eigenvalues (see appendix A) and the corresponding eigenstates
one obtains for
\begin{eqnarray}
&&V_1=\left|111\right\rangle+\left|1\overline{1}\overline{1}\right\rangle+\left|\overline{1}1\overline{1}\right\rangle
+\left|\overline{1}\overline{1}1\right\rangle\\
&&V_2=\left|\overline{1}\overline{1}\overline{1}\right\rangle+\left|\overline{1}11\right\rangle+
\left|1\overline{1}1\right\rangle+\left|11\overline{1}\right\rangle\\
&&T^{(3)}\left(\theta\right)\left(V_1,V_2\right)=e^{3m_{11}^{(+)}\theta}\left(V_1,V_2\right).
\end{eqnarray}
Along with $\left|222\right\rangle$ these assure
\begin{eqnarray}
&&\hbox{tr}\left(T^{(3)}\left(\theta\right)\right)=2e^{3m_{11}^{(+)}\theta}+1.
\end{eqnarray}

For $r=4$, the two corresponding combinations are
\begin{eqnarray}
&&V_1=\left|1111\right\rangle+\left|11\overline{1}\overline{1}\right\rangle+\left|1\overline{1}\overline{1}1\right\rangle
+\left|1\overline{1}1\overline{1}\right\rangle+\left|\overline{1}\overline{1}\overline{1}\overline{1}\right\rangle+
\left|\overline{1}\overline{1}11\right\rangle+\left|\overline{1}1\overline{1}1\right\rangle+
\left|\overline{1}11\overline{1}\right\rangle\nonumber\\
&&V_2=\left|111\overline{1}\right\rangle+\left|11\overline{1}1\right\rangle+\left|1\overline{1}11\right\rangle
+\left|\overline{1}111\right\rangle+\left|\overline{1}\overline{1}\overline{1}1\right\rangle+
\left|\overline{1}\overline{1}1\overline{1}\right\rangle+\left|\overline{1}1\overline{1}\overline{1}\right\rangle+
\left|1\overline{1}\overline{1}\overline{1}\right\rangle
\end{eqnarray}
with
\begin{equation}
T^{(4)}\left(\theta\right)\left(V_1,V_2\right)=e^{4m_{11}^{(+)}\theta}\left(V_1,V_2\right).
\end{equation}
Along with $\left|2222\right\rangle$ these assure
\begin{eqnarray}
&&\hbox{tr}\left(T^{(4)}\left(\theta\right)\right)=2e^{4m_{11}^{(+)}\theta}+1.
\end{eqnarray}

The general pattern is now visible. Indeed the general result is that (with relative coefficients in the sums below
being all unity as in the example above)
\begin{eqnarray}
&&V_e=\left(\hbox{sum of states with {\it even} number of} \left|1\right\rangle\right)\nonumber\\
&&V_0=\left(\hbox{sum of states with {\it odd} number of} \left|1\right\rangle\right)
\end{eqnarray}
give
\begin{equation}
T^{(r)}\left(\theta\right)\left(V_e,V_0\right)=e^{rm_{11}^{(+)}\theta}\left(V_e,V_0\right).
\end{equation}
Along with $\left|22\ldots 2\right\rangle$ they assure
\begin{eqnarray}
&&\hbox{tr}\left(T^{(r)}\left(\theta\right)\right)=2e^{rm_{11}^{(+)}\theta}+1.
\end{eqnarray}
All the reaming $\left(3^r-3\right)$ states are grouped unto subsets giving "zero-trace" multiplets of eigenvalues
involving roots of unity.

This a as far as we propose to go in explicitly constructing the eigenstates. We repeat that, as mentioned in appendix
A, we could obtain the complete sets for $r=1,2,3,4$ the $3,3^2,3^3,3^4$ eigenstates exhaustively.

Let us now compare our approach with that via algebraic Bethe ansatz (Refs. \cite{3,5,6} provide a
considerable number of references). For ready comparison we recapitulate the essential results for the relatively
simple and well studied case of $6$-vertex models. We follow the notation of ref. \cite{3} for the {\it ferroelectric
regime} in particular. Starting with the $4\times 4$ 6-vertex braid matrix and denoting the $N$-th order transfer
matrix blocks as ($N$, $r$ having {\it different} significance here as compared to our notations)
\begin{equation}
T^{(N)}\left(\theta\right)=\left(\begin{array}{cc}
  A\left(\theta\right) & B\left(\theta\right) \\
  C\left(\theta\right) & D\left(\theta\right) \\
\end{array}\right)
\end{equation}
with
\begin{equation}
\hbox{tr}\left(T^{(N)}\left(\theta\right)\right)=A\left(\theta\right)+D\left(\theta\right).
\end{equation}
The eigenvalues of this trace are extracted from the ansatz
\begin{equation}
\psi\left(\theta_1,\theta_2,\ldots,\theta_r\right)=B\left(\theta_1\right)B\left(\theta_2\right)
\cdots B\left(\theta_r\right)\left(\left|{1\atop 0}\right\rangle_1\otimes\left|{1\atop 0}\right\rangle_2\otimes\cdots\otimes
\left|{1\atop 0}\right\rangle_N\right).
\end{equation}
For the chosen regime, denoting
\begin{equation}
\lambda_j=i\left(\theta_j+\frac{\gamma}2\right)
\end{equation}
(where $\gamma$ is the single free parameter of $\widehat{R}$) the constraints on the parameters $\left(\theta_1,\ldots,\theta_r\right)$ for (3.31) to give eigenstates reduce (due
to the $Rtt$ algebra) to
\begin{equation}
\left[\frac{\sin\left(\lambda_j+i\gamma/2\right)}{\sin\left(\lambda_j-i\gamma/2\right)}\right]^N
=-\prod_{k=1}^r\frac{\sin\left(\lambda_j-\lambda_k+i\gamma\right)}{\sin\left(\lambda_j-\lambda_k-i\gamma\right)}.
\end{equation}
One has to find the solutions (in general complex) for $\left(\theta_1,\ldots,\theta_r\right)$ from these set of
nonlinear constraints.

For our case the invariance of the subspaces $S\left(r,k\right)$ defined in (2.10) under the action of $T^{(r)}
\left(\theta\right)$ (our $r$ being $N$ of (3.29)) clearly indicates the choice of
\begin{equation}
\left|22\ldots 2\right\rangle=\left|\begin{array}{c}
  0 \\
  1 \\
  0 \\
\end{array}\right\rangle\otimes\left|\begin{array}{c}
  0 \\
  1 \\
  0 \\
\end{array}\right\rangle\otimes\cdots\otimes\left|\begin{array}{c}
  0 \\
  1 \\
  0 \\
\end{array}\right\rangle
\end{equation}
with eigenvalue 1 (for all successive orders $r$) as the starting point. This is the subspace $S\left(r,r\right)$
with only one state.

From (3.1) to (3.28) we have defined and implemented operators which  acting on (3.34) moves stepwise through the subspaces
\begin{equation}
S\left(r,r\right),\,S\left(r,r-1\right),\ldots,\,S\left(r,1\right),\,S\left(r,0\right).
\end{equation}
For our class of higher dimensional structures, already for $N=3$, the state-labels $\left(1,2,\overline{1}\right)$
necessitate different types of actions on the index 2. Instead of a single complex $\theta$-dependent matrix
(like $B\left(\theta\right)$ of (3.31)) we have chosen and systematically implemented the operators $\left(U,A,D\right)$
to move through the sequence (3.35). Al each step our formalism leads to relatively simple {\it linear} constraints, (3.6)
giving the complete results for $S\left(r,r-1\right)$ being the simplest example. The results (3.26), (3.27) and (3.28)
give the complete trace. The results in appendix A (for $r=
1,2,3,4$ for spaces of $3,9,27,81$ dimensions respectively) give a fair idea of the structure of the eigenvalues.

The fact that one finally solves only sets of linear equations with simple constant coefficients is not evident directly
from (3.13), (3.14) and (3.15) for example. But the fact that one ends up only with eigenvalues of the form $e^{\mu\theta}$
(where, as in appendix A, $\mu$ is a linear function of the parameters $m_{ij}^{(\pm)}$) leads finally to such
constraints. For the $3^r-3$ eigenvalues of zero total trace one searches for multiplets formed by roots unity and hence
summing to zero. This also is very helpful in constructing eigenstates. The eigenvalues, for each $r$, are known to a
certain extent (through certainly not entirely) beforehand.

\section{Hamiltonians and conserved quantities ($N=3$)}
\setcounter{equation}{0}

We study here the role of our parameters in the sequence of conserved quantities, the first one in the sequence being
chain Hamiltonian \cite{3,4}. (Many source are cited in ref. \cite{3}.) Define
\begin{equation}
H_n=\left.\frac{\partial^n}{\partial \theta^n}\ln T^{(r)}\left(\theta\right)\right|_{\theta=0}.
\end{equation}
The commutativity of the transfer matrices $T\left(\theta\right)$, $T\left(\theta'\right)$ implies
\begin{equation}
\left[H_n,H_m\right]=0.
\end{equation}
If $H_1$ is regarded as the Hamiltonian of the system, there is infinite set of conserved quantities.

Using standard results \cite{3,4} and taking account of our normalization and the regularity, i.e.
\begin{equation}
\widehat{R}\left(0\right)=PR\left(0\right)=I,
\end{equation}
one obtains (since $\left(P\widehat{R}\left(0\right)\right)^{-1}\left(P\partial_{\theta}\widehat{R}
\left(\theta\right)\right)_{\theta=0}=\left(\partial_{\theta}\widehat{R}\left(\theta\right)\right)_{\theta=0}$)
\begin{equation}
H_1=\left(T^{(r)}\left(0\right)\right)^{-1}\left.\frac{\partial}{\partial \theta}T^{(r)}\left(\theta\right)
\right|_{\theta=0}=\sum_{k=1}^r I\otimes\cdots\otimes\widehat{R}\left(0\right)_{k,k+1}\otimes\cdots\otimes I.
\end{equation}
Note that due to the trace (circular boundary) constraint $k+1=r+1\approx 1$. Indeed starting with $r=2$, evaluating
directly and explicitly
\begin{equation}
\left(T^{(2)}\left(0\right)\right)^{-1}\left.\frac{\partial}{\partial\theta}T^{(2)}\left(\theta\right)\right|_{\theta=0}
\end{equation}
and setting
\begin{equation}
x_{\pm}=\frac 12\left(m_{11}^{(+)}\pm m_{11}^{(-)}\right),\qquad
y_{\pm}=\frac 12\left(m_{12}^{(+)}\pm m_{12}^{(-)}\right),\qquad
z_{\pm}=\frac 12\left(m_{21}^{(+)}\pm m_{21}^{(-)}\right),
\end{equation}
one obtains (writing $H$ for $H_1$ when $r=2$)
\begin{eqnarray}
&&H=\left(\begin{array}{ccccccccc}
  2x_+ & 0 & 0 & 0 & 0 & 0 & 0 & 0 & 2x_- \\
  0 & y_++z_+ & 0 & 0 & 0 & 0 & 0 & y_-+z_- & 0 \\
  0 & 0 & 2x_+ & 0 & 0 & 0 & 2x_- & 0 & 0 \\
  0 & 0 & 0 & y_++z_+ & 0 & y_-+z_- & 0 & 0 & 0 \\
  0 & 0 & 0 & 0 & 0 & 0 & 0 & 0 & 0 \\
  0 & 0 & 0 & y_-+z_- & 0 & y_++z_+ & 0 & 0 & 0 \\
  0 & 0 & 2x_- & 0 & 0 & 0 & 2x_+ & 0 & 0 \\
  0 & y_-+z_- & 0 & 0 & 0 & 0 & 0 & y_++z_+ & 0 \\
  2x_- & 0 & 0 & 0 & 0 & 0 & 0 & 0 & 2x_+ \\
\end{array}\right)\nonumber\\
&&\phantom{H}=\dot{\widehat{R}}\left(0\right)+P\dot{\widehat{R}}\left(0\right)P\nonumber\\
&&\phantom{H}=\dot{\widehat{R}}\left(0\right)_{12}+\dot{\widehat{R}}\left(0\right)_{21},
\end{eqnarray}
where, in evident notations,
\begin{eqnarray}
&&\dot{\widehat{R}}\left(0\right)_{12}=\left(x_+,\,y_+,\,x_+,\,z_+,\,0,\,z_+,\,x_+,\,y_+,\,x_+\right)_{\hbox{diag.}}+
\nonumber\\&&\phantom{\dot{\widehat{R}}\left(0\right)_{12}=}\left(x_-,\,y_-,\,x_-,\,z_-,\,0,\,z_-,\,x_-,\,y_-,
\,x_-\right)_{\hbox{anti-diag.}}
\end{eqnarray}
and
\begin{eqnarray}
&&\dot{\widehat{R}}\left(0\right)_{21}=\left(x_+,\,z_+,\,x_+,\,y_+,\,0,\,y_+,\,x_+,\,z_+,\,x_+\right)_{\hbox{diag.}}+
\nonumber\\&&\phantom{\dot{\widehat{R}}\left(0\right)_{12}=}\left(x_-,\,z_-,\,x_-,\,y_-,\,0,\,y_-,\,x_-,\,z_-,
\,x_-\right)_{\hbox{anti-diag.}}.
\end{eqnarray}
The expressions for (12) and (21) are related through the interchanges
\begin{equation}
\left(y_{\pm},z_{\pm}\right)\rightarrow \left(z_{\pm},y_{\pm}\right).
\end{equation}
The appearance of (21) in (4.7) is consistent with the remark below (4.4).

For higher derivatives one has
\begin{eqnarray}
&&\left.\frac{d^l}{d\theta^l}\widehat{R}\left(\theta\right)\right|_{\theta=0}=
\left(x_+^l,\,y_+^l,\,x_+^l,\,z_+^l,\,0,\,z_+^l,\,x_+^l,\,y_+^l,\,x_+^l\right)_{\hbox{diag.}}+\nonumber\\
&&\phantom{\left.\frac{d^l}{d\theta^l}\widehat{R}\left(\theta\right)\right|_{\theta=0}=}
\left(x_-^l,\,y_-^l,\,x_-^l,\,z_-^l,\,0,\,z_-^l,\,x_-^l,\,y_-^l,\,x_-^l\right)_{\hbox{anti-diag.}}.
\end{eqnarray}
For $H_2$ one now obtains, as compared to (4.4),
\begin{eqnarray}
&&H_2=\sum_{j\neq k} I\otimes\cdots\otimes I\otimes
\dot{\widehat{R}}\left(0\right)_{j,j+1}\otimes I\cdots\otimes I\otimes
\dot{\widehat{R}}\left(0\right)_{k,k+1}\otimes I\cdots\otimes I+\nonumber\\
&&\phantom{H_2=}\sum_{k}I\cdots\otimes I\otimes
\ddot{\widehat{R}}\left(0\right)_{k,k+1}\otimes I\cdots\otimes I.
\end{eqnarray}
Generalization to higher orders are carried out in evident fashion.

In section 5 of ref. \cite{1} in constructing $\theta$-expansions the $H$ defined in (5.1) is {\it precisely}
$\dot{\widehat{R}}\left(0\right)$ of (4.8) above generalized to all odd $N$, namely $N=3,\,5,\,7,\ldots$. There it was
noted (eq. (5.9) of ref. \cite{1}),
\begin{equation}
\left[H_{12}+H_{23},\left[H_{12},H_{23}\right]\right]=0,
\end{equation}
where $H_{12}=H\otimes I$ and $H_{23}=I\otimes H$. This vanishing double commutator is the simplest version of the
Reshetikhin condition given in eqs. (3.20)of ref. \cite{2} as
\begin{equation}
\left[H_{12}+H_{23},\left[H_{12},H_{23}\right]\right]=X_{12}-X_{23},
\end{equation}
the r.h.s. being the difference of two two-point-quantities. In (4.13) the r.h.s. is simply zero.

\section{Potential for factorizable $S$-matrices and Cayley transforms ($N=3$)}
\setcounter{equation}{0}

Potentials for scattering of bosons or fermions with quadratic interaction terms (sec. 3 of ref. 2 and sec. 1 of ref. 3
provide more references) can correspond to factorizable $S$-matrices (factorizable into two particle scatterings,
independently of the chosen order of the latter ones) provided that such potentials are inverse Cayley
transforms of Yang-Baxter matrices of appropriate dimensions, i.e. $V$ being the potential (for a chosen helicity fixing
the sign of $\theta$)
\begin{equation}
-iV=\left(R\left(\theta\right)-\lambda\left(\theta\right)I\right)^{-1}
\left(R\left(\theta\right)+\lambda\left(\theta\right)I\right).
\end{equation}
As compared to refs. 2, 3 we display explicitly a free normalization factor
\begin{equation}
\left(\lambda\left(\theta\right)\right)^{-1}R\left(\theta\right).
\end{equation}
Our multiparametric case shows clearly that though the normalization (if well-defined) trivially cancels in the YB or
the braid equation it must be compatible with the existence of the inverse of
\begin{equation}
\left(\lambda^{-1}\left(\theta\right)R\left(\theta\right)-I\right).
\end{equation}
We will find that
\begin{equation}
\lambda\left(\theta\right)\neq
\left(1,e^{m_{11}^{(\pm)}\theta},\pm e^{\frac 12(
m_{12}^{(\pm)}+m_{21}^{(\pm)})\theta}\right).
\end{equation}
The inverse (when $\widehat{R}\left(\theta\right)$ is given by (1.10))
\begin{equation}
\left(\widehat{R}\left(\theta\right)-\lambda'\left(\theta\right)I\right)^{-1}
\end{equation}
can be shown to exist for
\begin{equation}
\lambda'\left(\theta\right)\neq
\left(1,e^{m_{11}^{(\pm)}\theta},e^{m_{12}^{(\pm)}\theta},e^{m_{21}^{(\pm)}\theta}\right).
\end{equation}
The significance of (5.6) is simple, the r.h.s. exhibiting simply the coefficient of the projector
\begin{equation}
\left(P_{22},P_{11}^{(\pm)},P_{12}^{(\pm)},P_{21}^{(\pm)}\right)
\end{equation}
in (1.10).

Diagonalizing $\widehat{R}\left(\theta\right)$ the situation becomes particularly transparent (eqs. (3.4), (3.5) of
ref. 1). $M$ being given by (3.5) of ref. 1,
\begin{eqnarray}
&&M\left(\widehat{R}\left(\theta\right)-\lambda'\left(\theta\right)I\right)M^{-1}=\nonumber\\
&&\left(e^{m_{11}^{(+)}\theta},e^{m_{12}^{(+)}\theta},e^{m_{11}^{(+)}\theta},
e^{m_{21}^{(+)}\theta},1,e^{m_{21}^{(-)}\theta},e^{m_{11}^{(-)}\theta},e^{m_{12}^{(-)}\theta},
e^{m_{11}^{(-)}\theta}\right)_{\hbox{diag.}}-\lambda'\left(\theta\right)I.
\end{eqnarray}
When $\lambda'\left(\theta\right)$ is equal to any one of the eigenvalues (including 1) the determinant of
$\left(\widehat{R}\left(\theta\right)-\right.$ $\left.\lambda'\left(\theta\right)I\right)$ vanishes. Hence (5.6).

For (5.1) one requires invertibility of
\begin{equation}
P\widehat{R}\left(\theta\right)-\lambda\left(\theta\right)I.
\end{equation}
The action of $P$ finally leads to (5.4) rather than (5.6).

Defining $X$ through
\begin{equation}
\left(R\left(\theta\right)-\lambda\left(\theta\right)I\right)X=I,
\end{equation}
we present below the explicit form of $X$ for our $N=3$ case.
\begin{equation}
X=\left(\begin{array}{ccccccccc}
  x_1 & 0 & 0 & 0 & 0 & 0 & 0 & 0 & x_8 \\
  0 & x_2 & 0 & x_6 & 0 & x_7 & 0 & x_4 & 0 \\
  0 & 0 & x_3 & 0 & 0 & 0 & x_9 & 0 & 0 \\
  0 & x_{10} & 0 & x_2 & 0 & x_4 & 0 & x_{11} & 0 \\
  0 & 0 & 0 & 0 & x_5 & 0 & 0 & 0 & 0 \\
  0 & x_{11} & 0 & x_4 & 0 & x_2 & 0 & x_{10} & 0 \\
  0 & 0 & x_9 & 0 & 0 & 0 & x_3 & 0 & 0 \\
  0 & x_4 & 0 & x_7 & 0 & x_6 & 0 & x_2 & 0 \\
  x_8 & 0 & 0 & 0 & 0 & 0 & 0 & 0 & x_1 \\
\end{array}\right),
\end{equation}
where (writing $\lambda$ for $\lambda\left(\theta\right)$)
\begin{eqnarray}
&&x_{1}=\frac 12\left(\frac
1{e^{m_{11}^{(+)}\theta}-\lambda}+\frac
1{e^{m_{11}^{(-)}\theta}-\lambda}\right),\qquad x_8=\frac 12\left(\frac 1{e^{m_{11}^{(+)}\theta}-\lambda}-\frac
1{e^{m_{11}^{(-)}\theta}-\lambda}\right),\nonumber\\
&&x_2=\frac \lambda 2\left(\frac
1{e^{(m_{12}^{(+)}+m_{21}^{(+)})\theta}-\lambda^2}+\frac
1{e^{(m_{12}^{(-)}+m_{21}^{(-)})\theta}-\lambda^2}\right),\nonumber\\
&&x_4=\frac \lambda 2\left(\frac
1{e^{(m_{12}^{(+)}+m_{21}^{(+)})\theta}-\lambda^2}-\frac
1{e^{(m_{12}^{(-)}+m_{21}^{(-)})\theta}-\lambda^2}\right),\nonumber\\
&&x_{3}=\frac 12\left(\frac
1{e^{m_{11}^{(+)}\theta}-\lambda}-\frac
1{e^{m_{11}^{(-)}\theta}+\lambda}\right),\qquad x_9=\frac 12\left(\frac 1{e^{m_{11}^{(+)}\theta}-\lambda}+\frac
1{e^{m_{11}^{(-)}\theta}+\lambda}\right),\nonumber\\
&&x_5=\frac 1{1-\lambda},\nonumber\\
&&x_6=\frac 12\left(\frac{e^{m_{21}^{(+)}\theta}}{e^{(m_{12}^{(+)}+m_{21}^{(+)})\theta}-\lambda^2}+\frac
{e^{m_{21}^{(-)}\theta}}{e^{(m_{12}^{(-)}+m_{21}^{(-)})\theta}-\lambda^2}\right),\nonumber\\
&&x_7=\frac 12\left(\frac
{e^{m_{21}^{(+)}\theta}}{e^{(m_{12}^{(+)}+m_{21}^{(+)})\theta}-\lambda^2}-\frac
{e^{m_{21}^{(-)}\theta}}{e^{(m_{12}^{(-)}+m_{21}^{(-)})\theta}-\lambda^2}\right),\nonumber\\
&&x_{10}=\frac 12\left(\frac{e^{m_{12}^{(+)}\theta}}{e^{(m_{12}^{(+)}+m_{21}^{(+)})\theta}-\lambda^2}+\frac
{e^{m_{12}^{(-)}\theta}}{e^{(m_{12}^{(-)}+m_{21}^{(-)})\theta}-\lambda^2}\right),\nonumber\\
&&x_{11}=\frac 12\left(\frac{e^{m_{12}^{(+)}\theta}}{e^{(m_{12}^{(+)}+m_{21}^{(+)})\theta}-\lambda^2}-\frac
{e^{m_{12}^{(-)}\theta}}{e^{(m_{12}^{(-)}+m_{21}^{(-)})\theta}-\lambda^2}\right).
\end{eqnarray}
Now, from (5.1), (5.10) and (5.11),
\begin{eqnarray}
&&-iV=\left(R\left(\theta\right)-\lambda\left(\theta\right)I\right)^{-1}
\left(R\left(\theta\right)+\lambda\left(\theta\right)I\right)\nonumber\\
&&\phantom{-iV}=\left(R\left(\theta\right)-\lambda\left(\theta\right)I\right)^{-1}
\left(R\left(\theta\right)-\lambda\left(\theta\right)I+2\lambda\left(\theta\right)I\right)\nonumber\\
&&\phantom{-iV}=X\left(X^{-1}+2\lambda\left(\theta\right)I\right)\nonumber\\
&&-iV=I+2\lambda\left(\theta\right)X.
\end{eqnarray}
From (5.12) it is evident that $X$ is well-defined only when (5.4) is satisfied. We have thus obtained explicitly,
for $N=3$, the potential leading to a factorizable $S$-matrix. The role played by our parameters is now displayed.

Note that any $\lambda\left(\theta\right)$ satisfying (5.4) can be implemented. One many choose to display this
dependence on $\lambda$ by denoting the potential as $V\left(\lambda\right)$. With our $V\left(\lambda\right)$ one now
considers the fermionic lagrangian
\begin{equation}
{\cal L}=\int dx\left(i\overline{\psi}_a\gamma_\nu\partial_\nu\psi_a-g\left(\overline{\psi}_a\gamma_\nu\psi_c\right)
V_{ab,cd}\left(\overline{\psi}_b\gamma_\nu\partial_\nu\psi_d\right)\right),
\end{equation}
where
\begin{equation}
V=\sum_{ab,cd}\left(V_{ab,cd}\right)\left(ab\right)\otimes\left(cd\right).
\end{equation}
There is an analogous, simpler,formulation for bosons. We will not further analyze the consequences of our $V$. But it
should be compared to the detailed studies of the solutions obtained in refs. \cite{7,8}.

\section{$N>3$}
\setcounter{equation}{0}

So far we have studied the case $N=3$ in detail. Now we indicate
briefly the crucial new features arising for $N>3$. Many aspects
are conserved also, as will be pointed out.

The first major feature is the generalization of (2.1). For $N=2p-1$, one obtains
\begin{equation}
\hbox{tr}\left(T^{(r)}\left(\theta\right)\right)=2\left(e^{rm_{11}^{(+)}\theta}+e^{rm_{22}^{(+)}\theta}+\cdots+
e^{rm_{p-1,p-1}^{(+)}\theta}\right)+1,\qquad p=2,\,3,\,4,\ldots
\end{equation}
There are
\begin{equation}
2\left(p-1\right)+1=2p-1=N
\end{equation}
terms. An explanation, promised below (2.1), is as follows. {\it Only the diagonal blocks $t_{ii}\left(\theta\right)$ have
diagonal terms}. Thus in (2.5) only $t_{11}\left(\theta\right)$, $t_{22}\left(\theta\right)$, $t_{\overline{1}
\overline{1}}\left(\theta\right)$, has non zero elements on the diagonal, their sum being
\begin{equation}
\hbox{tr}\left(T^{(1)}\left(\theta\right)\right)=2\left(a_++a_-\right)+1=2 e^{m_{11}^{(+)}\theta}+1,\qquad (r=1)
\end{equation}
For $r=2$ (and $N=3$) one obtains from the coproduct structure
\begin{eqnarray}
&&\hbox{tr}\left(T^{(2)}\left(\theta\right)\right)=\hbox{tr}\left(\left(a_+a_+,0,a_+a_-\right)_{\hbox{diag.}}
+\left(0,1,0\right)_{\hbox{diag.}}+\left(a_-a_-,0,a_-a_+\right)_{\hbox{diag.}}+\right.\nonumber\\
&&\phantom{\hbox{tr}\left(T^{(2)}\left(\theta\right)\right)=}\left.\left(a_-a_+,0,a_-a_-\right)_{\hbox{diag.}}
+\left(a_+a_-,0,a_+a_+\right)_{\hbox{diag.}}\right)+\nonumber\\
&&\phantom{\hbox{tr}\left(T^{(2)}\left(\theta\right)\right)=}\hbox{tr}\left(\hbox{blocks with nondiagonal terms only}\right)
\nonumber\\
&&\phantom{\hbox{tr}\left(T^{(2)}\left(\theta\right)\right)}=2\left(a_++a_-\right)^2+1\nonumber\\
&&\hbox{tr}\left(T^{(2)}\left(\theta\right)\right)=2e^{2m_{11}^{(+)}\theta}+1
\end{eqnarray}
and so on.

For $n>3$ the basic features persist. Along with the crucial constraint (1.8)
\begin{equation}
m_{ij}^{(\epsilon)}=m_{i\overline{j}}^{(\epsilon)},\qquad (\overline{j}=2p-j)
\end{equation}
which symmetrizes the blocks on the diagonal (generalizing (1.10)), the final result is
\begin{eqnarray}
&&\hbox{tr}\left(T^{(r)}\left(\theta\right)\right)=2\sum_{i=1}^{p-1}\left(a_{ii}^{(+)}+a_{ii}^{(-)}\right)^r+1\\
&&\phantom{\hbox{tr}\left(T^{(r)}\left(\theta\right)\right)}=2\sum_{i=1}^{p-1}e^{m_{ii}^{(+)}\theta}+1.
\end{eqnarray}
The result is obtained directly by looking closely at the structure of the matrices concerned, {\it without constructing
eigenstates and their eigenvalues}. We have checked (6.7) directly and explicitly for arbitrary $p$ (appendix A).
But this result has profound consequence on the spectrum of the eigenvalues for each $r$. Given $\left(N,r\right)$ and
the coproduct rule,
\begin{enumerate}
    \item the number of eigenvalues for $T^{(r)}\left(\theta\right)=N^r$
    \item the number of eigenvalues contributing in the trace=$N$
    \item the remaining $\left(N^r-N\right)$ eigenvalues must sum to give zero contribution in the trace.
\end{enumerate}
For $N=3$ we have shown (appendix A) how this constraint is satisfied via the multiplet structures
\begin{eqnarray}
&&e^{\mu\theta}\left(1,\omega_{(l)},\omega_{(l)}^2,\ldots,\omega_{(l)}^{l-1}\right),\qquad \omega_l=e^{\frac{2\pi i}l}\nonumber\\
&&1+\omega_{(l)}+\omega_{(l)}^2+\cdots+\omega_{(l)}^{l-1}=0,
\end{eqnarray}
where $\mu$ is linear in $m_{ij}^{(\pm)}$ and for
\begin{eqnarray}
&&r=2,\qquad l=2\nonumber\\
&&r=3,\qquad l=3\nonumber\\
&&r=4,\qquad l=2,\,4
\end{eqnarray}
and so on.

This multiplet structure involving roots of unity can also be shown to be carried over for $N>3$ explicitly. But apart
from the fact that the number of eigenstates and the number of states in the linear combinations giving eigenstates
increase very fast there are no other basic difficulties.

For $N=5$, for example, generalizing the basis (2.7), for $p=3$, to
\begin{equation}
\left(\left|\begin{array}{c}
  1 \\
  0 \\
  0 \\
  0 \\
  0 \\
\end{array}\right\rangle,\left|\begin{array}{c}
  0 \\
  1 \\
  0 \\
  0 \\
  0 \\
\end{array}\right\rangle,\left|\begin{array}{c}
  0 \\
  0 \\
  1 \\
  0 \\
  0 \\
\end{array}\right\rangle,\left|\begin{array}{c}
  0 \\
  0 \\
  0 \\
  1 \\
  0 \\
\end{array}\right\rangle,\left|{\begin{array}{c}
  0 \\
  0 \\
  0 \\
  0 \\
  1 \\
\end{array}}\right\rangle\right)\equiv\left(\left|1\right\rangle,\left|2\right\rangle,\left|3\right\rangle,
\left|\overline{2}\right\rangle,\left|\overline{1}\right\rangle,\right)
\end{equation}
one again obtains subspaces stable under the action of $T^{(r)}\left(\theta\right)$
\begin{equation}
S\left(r,k\right)\qquad \left(k=0,1,\ldots,r\right),
\end{equation}
where $k$ is now the multiplicity of the index 3. The operator structures (3.1), (3.2), (3.3) are now generalized, in
terms of operators $\left(t_{31}\left(\theta\right),t_{3\overline{1}}\left(\theta\right)\right)$, $\left(t_{32}
\left(\theta\right),t_{3\overline{2}}\left(\theta\right)\right)$, $t_{33}\left(\theta\right)$ to construct
\begin{eqnarray}
&&U_1\left|3\right\rangle=\left|1\right\rangle,\qquad U_2\left|3\right\rangle=\left|2\right\rangle,\nonumber\\
&&D_1\left|3\right\rangle=\left|\overline{1}\right\rangle,\qquad D_2\left|3\right\rangle=\left|\overline{2}\right
\rangle,\nonumber\\
&&A\left|3\right\rangle=\left|3\right\rangle.
\end{eqnarray}
The subspace $S\left(r,r\right)$ is still given by a single state
\begin{equation}
T^{(r)}\left(\theta\right)\left|33\ldots 3\right\rangle=\left|33\ldots 3\right\rangle.
\end{equation}
The subspace $S\left(r,r-1\right)$ is now spanned by $4r$ (instead of $2r$ for $N=3$) states and is easily diagonalized.

The stepwise generalization for $N=7,\,9,\ldots$ is now fairly evident. The dimensions of $S\left(r,k\right)$ is given
by a generalization of (2.11) by the successive coefficients in
\begin{equation}
\left(x+2\left(p-1\right)\right)^r=1\cdot x^r+2r\left(p-1\right)x^{r-1}+\cdots+
\left(2\left(p-1\right)\right)^{r-k}\binom{r}{r-k}x^k+\cdots+\left(2\left(p-1\right)\right)^{r}.
\end{equation}
For $x=1$, one gets the total dimension
\begin{equation}
\left(2p-1\right)^r=N^{r}.
\end{equation}
The generalization of (3.26) is also fairly direct. When the order $r$ of $T^{(r)}\left(\theta\right)$ is a prime number
there is an amusing encounter with a theorem of Fermat in considering our multiplet structures. This is discussed in
appendix B.

The generalization of the structure of the Hamiltonian of sec. 4 involving $\dot{\widehat{R}}\left(0\right)$ of (4.8)
is particulary straightforward. But even for $N=5$ we have 24 non-zero terms on the diagonal and as many on the
antidiagonal.

The potential (5.1) for $N=5,\,7,\ldots$ now involve inversions of $N^2$ dimensional matrices $\left(R\left(\theta\right)
-\lambda\left(\theta\right)I\right)$. This is again straightforward, given our specific structure of $\widehat{R}
\left(\theta\right)$, but evidently lengthly.

Apart from such general indications as presented above systematic studies for cases $N>3$ are beyond the scope of this
work. In particular possible substructures in each subspace $S\left(r,k\right)$ corresponding to multiplicities of
different indices (say $(1,\overline{1})$ and $(2,\overline{2})$ for $N=5$) should be formulated with care.

\section{Generalization of the nested sequence of projectors}
\setcounter{equation}{0}

The sequence of projectors (1.3), forming a complete ortho-normalized basis admits the more general parametrization
displayed below (for odd $N$)
\begin{eqnarray}
&& P_{pp}=(pp)\otimes (pp),\nonumber\\
&&\left(u_{pi}+u_{pi}^{-1}\right)P_{pi\left(\pm\right)}=(pp)\otimes\left[u_{pi}^{\pm 1}(ii)+
u_{pi}^{\mp 1}(\overline{i}\overline{i})\pm \left(v_{pi}(i\overline{i})+
v_{pi}^{-1}(\overline{i}i)\right)\right],\nonumber\\
&&\left(u_{ip}+u_{ip}^{-1}\right)P_{ip\left(\pm\right)}=\left[u_{ip}^{\pm 1}(ii)+
u_{ip}^{\mp 1}(\overline{i}\overline{i})\pm \left(v_{ip}(i\overline{i})+
v_{ip}^{-1}(\overline{i}i)\right)\right]\otimes(pp),\\
&&\left(u_{ij}+u_{ij}^{-1}\right)P_{ij\left(\epsilon\right)}=u_{ij}^{\pm 1}(ii)\otimes(jj)+
u_{ij}^{\mp 1}(\overline{i}\overline{i})\otimes(\overline{j}\overline{j})+
\pm\left[v_{ij}(i\overline{i})\otimes(j\overline{j})+v_{ij}^{-1}(\overline{i}i)\otimes(\overline{j}j)\right],\nonumber\\
&&\left(u_{i\overline{j}}+u_{i\overline{j}}^{-1}\right)P_{i\overline{j}\left(\epsilon\right)}=u_{i\overline{j}}^{\pm 1}
(ii)\otimes(\overline{j}\overline{j})+u_{i\overline{j}}^{\mp 1}(\overline{i}\overline{i})\otimes(jj)+
\pm\left[v_{i\overline{j}}(i\overline{i})\otimes(\overline{j}j)+v_{i\overline{j}}^{-1}(\overline{i}i)\otimes
(j\overline{j})\right],\nonumber
\end{eqnarray}
where the supplementary parameters introduced are compatible with the orthonormality and completeness conditions (1.5).
For even $N$ also an analogues parametrization can be introduced. Thus $6$-vertex and $8$-vertex projector basis (given
in (6.1) of ref. 1) can be generalized to
\begin{eqnarray}
&&\left(u_{11}+u_{11}^{-1}\right)P_{11\left(\pm\right)}=
\left(\begin{array}{cccc}
  u_{11}^{\pm 1} & 0 & 0 & \pm v_{11} \\
  0 & 0 & 0 & 0 \\
  0 & 0 & 0 & 0 \\
  \pm v_{11}^{-1} & 0 & 0 & u_{11}^{\mp 1} \\
\end{array}\right)\nonumber\\
&&\left(u_{1\overline{1}}+u_{1\overline{1}}^{-1}\right)P_{1\overline{1}\left(\pm\right)}=
\left(\begin{array}{cccc}
  0 & 0 & 0 & 0 \\
  0 & u_{1\overline{1}}^{\pm 1} & \pm v_{1\overline{1}} & 0 \\
  0 & \pm v_{1\overline{1}}^{-1} & u_{1\overline{1}}^{\mp 1} & 0 \\
  0 & 0 & 0 & 0 \\
\end{array}\right)
\end{eqnarray}
Braid matrices on such bases of projectors and associated statistical models can be studied systematically. Such a study
will be presented elsewhere with suitable restrictions on parameters for specific solutions. Setting $u_{ab}=v_{ab}=1$ for
all values of the indices one recovers the projectors of (1.3).

\section{Discussion}
\setcounter{equation}{0}

Our results presented above, are limited to formal study of the transfer matrix, construction of chain Hamiltonians and
potentials corresponding to factorizable $S$-matrices. Adequate study of the, consequences of the features obtained, of
their deeper significance remains to be done. The role of our parameters should be analyzed in various domains for
comparison with corresponding features of well-known statistical models \cite{3,9,10,11}.

Certain features are, of course, immediately available for our case. Thus the free energy (defined with the opposite
sign in ref. \cite{10}), is given by the maximum eigenvalue $\left(\theta>0\right)$ as
\begin{equation}
f=-\lim_{r\rightarrow \infty}\frac 1r\ln e^{rm_{11}^{(+)}\theta}=-m_{11}^{(+)}\theta
\end{equation}
if we choose, say, the order
\begin{equation}
m_{11}^{(+)}>m_{22}^{(+)}>\cdots > m_{p-1,p-1}^{(+)}.
\end{equation}
In our case results depend on the sector of the parameters selected (their ordering). The second largest eigenvalue,
also of particular interest, is again directly obtained once the ordering is fixed.

Correlation functions are of major interest and a domain of intense activity \cite{12,13}. Here our model can have quite
interesting consequences. This aspect also remains to be explored.

We intend to continue our study elsewhere. But we consider the series of remarkable features presented here to be
sufficiently rich in content. They open up a significantly different domain, as compared to standard, well known cases.

\vskip 0.5cm

\noindent{\bf Acknowledgments:} {\em It is a pleasure to thank Jean Lascoux and Alain Lascoux for discussions concerning
the theme of appendix B. One of us (BA) wants to thank Patrick Mora for a Kind invitation at Ecole Polytechnique. He is
also very grateful to the members of the group for their warm hospitality. This work is supported by a grant of
"La Fondation Charles de Gaulle".}

\newpage

\begin{appendix}

\section{\LARGE Eigenvalues of $T^{(r)}\left(\theta\right)$ for $N=3$, $r=1,2,3,4$ and direct
construction of trace for $N>3$}
\setcounter{equation}{0}

For each $r$, the eigenvalues of $T^{(r)}\left(\theta\right)$ is given systematically for the subspaces
$S\left(r,k\right)$ defined in (2.10), for $k=0,1,\ldots,r$.

\paragraph{$\underline{r=1\,(\dim\,3)}$:}
\begin{eqnarray}
&&\phantom{S\left(1,0\right)\,\left(\dim\,2\right):\qquad}\hbox{eigenvalues}\nonumber\\
&&S\left(1,0\right)\,\left(\dim\,2\right):\qquad e^{m_{11}^{(+)}\theta}\left(1,1\right)\\
&&S\left(1,1\right)\,\left(\dim\,1\right):\qquad 1\\
&&\underline{\hbox{tr}\left(T^{(1)}\left(\theta\right)\right)=2e^{m_{11}^{(+)}\theta}+1}.
\end{eqnarray}

\paragraph{$\underline{r=2\, (\dim\, 9)}$:}
\begin{eqnarray}
&&S\left(2,0\right)\,\left(\dim\,4\right):\qquad e^{2m_{11}^{(+)}\theta}\left(1,1\right)\\
&&\phantom{S\left(2,0\right)\,\left(\dim\,4\right):\qquad} e^{2m_{11}^{(-)}\theta}\left(1,-1\right)\nonumber\\
&&S\left(2,1\right)\,\left(\dim\,4\right):\qquad e^{(m_{12}^{(+)}+m_{21}^{(+)})\theta}
\left(1,-1\right)\\
&&\phantom{S\left(2,0\right)\,\left(\dim\,4\right):\qquad} e^{(m_{12}^{(-)}+m_{21}^{(-)})\theta}
\left(1,-1\right)\nonumber\\
&&S\left(2,2\right)\,\left(\dim\,1\right):\qquad 1\\
&&\underline{\hbox{tr}\left(T^{(2)}\left(\theta\right)\right)=2e^{2m_{11}^{(+)}\theta}+1}.
\end{eqnarray}

\paragraph{$\underline{r=3\,(\dim\, 27)}$:}
\begin{eqnarray}
&&S\left(3,0\right)\,\left(\dim\,8\right):\qquad e^{3m_{11}^{(+)}\theta}\left(1,1\right)\nonumber\\
&&\phantom{S\left(3,0\right)\,\left(\dim\,8\right):\qquad} e^{(m_{11}^{(+)}+2m_{11}^{(-)})\theta}
\left(1,e^{\frac{2\pi i}3},e^{\frac{2\pi i}3\cdot 2}\right) \qquad \hbox{[2 times]}\\
&&S\left(3,1\right)\,\left(\dim\,12\right):\qquad \left(\begin{array}{c}
  e^{(m_{11}^{(+)}+m_{12}^{(+)}+m_{21}^{(+)})\theta}\\
  e^{(m_{11}^{(+)}+m_{12}^{(-)}+m_{21}^{(-)})\theta} \\
  e^{(m_{11}^{(-)}+m_{12}^{(+)}+m_{21}^{(-)})\theta} \\
  e^{(m_{11}^{(-)}+m_{12}^{(-)}+m_{21}^{(+)})\theta} \\
\end{array}\right)\left(1,e^{\frac{2\pi i}3},e^{\frac{2\pi i}3\cdot 2}\right)\\
&&S\left(3,2\right)\,\left(\dim\,6\right):\qquad\left(\begin{array}{c}
  e^{(m_{12}^{(+)}+m_{21}^{(+)})\theta} \\
  e^{(m_{12}^{(-)}+m_{21}^{(-)})\theta} \\
\end{array}\right)\left(1,e^{\frac{2\pi i}3},e^{\frac{2\pi i}3\cdot 2}\right)\\
&&S\left(3,3\right)\,\left(\dim\,1\right):\qquad 1\\
&&\underline{\hbox{tr}\left(T^{(3)}\left(\theta\right)\right)=2e^{3m_{11}^{(+)}\theta}+1}.
\end{eqnarray}

\paragraph{$\underline{r=4\, (\dim\, 81)}$:}
\begin{eqnarray}
&&S\left(4,0\right)\,\left(\dim\,16\right):\qquad e^{4m_{11}^{(+)}\theta}\left(1,1\right)\nonumber\\
&&\phantom{S\left(4,0\right)\,\left(\dim\,16\right):\qquad} e^{2(m_{11}^{(+)}+m_{11}^{(-)})\theta}
\left(1,e^{\frac{2\pi i}4},e^{\frac{2\pi i}4\cdot 2},e^{\frac{2\pi i}4\cdot 3}\right)\qquad\hbox{[3 times]}\nonumber\\
&&\phantom{S\left(4,0\right)\,\left(\dim\,16\right):\qquad} e^{4m_{11}^{(-)}\theta}\left(1,-1\right)\\
&&S\left(4,1\right)\,\left(\dim\,32\right):\qquad e^{(m_{11}^{(+)}+m_{11}^{(-)}+m_{12}^{(+)}+m_{21}^{(-)})
\theta}\left(1,e^{\frac{2\pi i}4},e^{\frac{2\pi i}4\cdot 2},e^{\frac{2\pi i}4\cdot 3}\right)\qquad\hbox{[2 times]}
\nonumber\\
&&\phantom{S\left(4,1\right)\,\left(\dim\,32\right):\qquad}e^{(m_{11}^{(+)}+m_{11}^{(-)}+m_{12}^{(-)}+m_{21}^{(+)})
\theta}\left(1,e^{\frac{2\pi i}4},e^{\frac{2\pi i}4\cdot 2},e^{\frac{2\pi i}4\cdot 3}\right)\qquad\hbox{[2 times]}
\nonumber\\
&&\phantom{S\left(4,1\right)\,\left(\dim\,32\right):\qquad}\left(\begin{array}{c}
  e^{(2m_{11}^{(+)}+m_{12}^{(+)}+m_{21}^{(+)})\theta} \\
  e^{(2m_{11}^{(+)}+m_{12}^{(-)}+m_{21}^{(-)})\theta} \\
  e^{(2m_{11}^{(-)}+m_{12}^{(+)}+m_{21}^{(+)})\theta} \\
  e^{(2m_{11}^{(-)}+m_{12}^{(-)}+m_{21}^{(+)})\theta} \\
\end{array}\right)
\left(1,e^{\frac{2\pi i}4},e^{\frac{2\pi i}4\cdot 2},e^{\frac{2\pi i}4\cdot 3}\right)\\
&&S\left(4,2\right)\,\left(\dim\,24\right):\qquad\left(\begin{array}{c}
  e^{2(m_{12}^{(+)}+m_{21}^{(+)})\theta} \\
  e^{2(m_{12}^{(-)}+m_{21}^{(-)})\theta} \\
\end{array}\right)\left(1,-1\right)\nonumber\\
&&\phantom{S\left(4,2\right)\,\left(\dim\,24\right):\qquad}
e^{(m_{12}^{(+)}+m_{12}^{(-)}+m_{21}^{(+)}+m_{21}^{(-)})\theta}\left(1,e^{\frac{2\pi i}4},e^{\frac{2\pi i}4\cdot 2},
e^{\frac{2\pi i}4\cdot 3}\right)\nonumber\\
&&\phantom{S\left(4,2\right)\,\left(\dim\,24\right):\qquad}\left(\begin{array}{c}
  e^{(m_{11}^{(+)}+m_{12}^{(+)}+m_{21}^{(+)})\theta} \\
  e^{(m_{11}^{(+)}+m_{12}^{(-)}+m_{21}^{(-)})\theta} \\
  e^{(m_{11}^{(-)}+m_{12}^{(+)}+m_{21}^{(-)})\theta} \\
  e^{(m_{11}^{(-)}+m_{12}^{(-)}+m_{21}^{(+)})\theta} \\
\end{array}\right)
\left(1,e^{\frac{2\pi i}4},e^{\frac{2\pi i}4\cdot 2},e^{\frac{2\pi i}4\cdot 3}\right)\\
&&S\left(4,3\right)\,\left(\dim\,8\right):\qquad\left(\begin{array}{c}
  e^{(m_{12}^{(+)}+m_{21}^{(+)})\theta} \\
  e^{(m_{12}^{(-)}+m_{21}^{(-)})\theta} \\
\end{array}\right)\left(1,e^{\frac{2\pi i}4},e^{\frac{2\pi i}4\cdot 2},e^{\frac{2\pi i}4\cdot 3}\right)\\
&&S\left(4,4\right)\,\left(\dim\,1\right):\qquad 1\\
&&\underline{\hbox{tr}\left(T^{(4)}\left(\theta\right)\right)=2e^{4m_{11}^{(+)}\theta}+1}.
\end{eqnarray}
Now we indicate briefly the direct construction of trace for all $N$, without constructing the full set of eigenvalues
explicitly.

Set, for the coefficients on the diagonal and the anti-diagonal respectively
\begin{equation}
d_{ij}=\frac 12\left(e^{m_{ij}^{(+)}\theta}+e^{m_{ij}^{(-)}\theta}\right)=d_{i\overline{j}},\qquad
a_{ij}=\frac 12\left(e^{m_{ij}^{(+)}\theta}-e^{m_{ij}^{(-)}\theta}\right)=a_{i\overline{j}},
\end{equation}
where
\begin{equation}
i=1,\,2,\ldots,\,p-1\qquad \overline{i}=(2p-1),\ldots,\,(p+1)\qquad p=\frac 12\left(N+1\right).
\end{equation}
From (1.7) and (1.15)
\begin{eqnarray}
&&t\left(\theta\right)=\sum_i\left((pi)\otimes\left(d_{ip}(ip)+a_{ip}(\overline{i}p)\right)+
(p\overline{i})\otimes\left(a_{ip}(ip)+d_{ip}(\overline{i}p)\right)+\right.\nonumber\\
&&\phantom{t\left(\theta\right)=}\left.(ip)\otimes\left(d_{pi}(pi)+a_{pi}(p\overline{i})\right)+(\overline{i}p)\otimes\left(a_{pi}(pi)+
d_{pi}(p\overline{i})\right)\right)\nonumber\\
&&\phantom{t\left(\theta\right)=}\sum_{i,j}\left((ji)\otimes\left(d_{ij}(ij)+a_{ij}(\overline{i}\overline{j})\right)+
(\overline{j}\overline{i})\otimes\left(a_{ij}(ij)+d_{ij}(\overline{i}\overline{j})\right)+\right.\nonumber\\
&&\phantom{t\left(\theta\right)=}\left.(j\overline{i})\otimes\left(d_{ij}(\overline{i}j)+a_{ij}(i\overline{j})\right)+
(\overline{j}i)\otimes\left(d_{ij}(\overline{i}j)+a_{ij}(i\overline{j})\right)\right)\nonumber\\
&&\phantom{t\left(\theta\right)=}\sum_{i}\left((i\overline{i})\otimes\left(d_{ii}(\overline{i}i)+a_{ii}
(i\overline{i})\right)+(\overline{i}i)\otimes\left(a_{ii}(\overline{i}i)+d_{ii}(i\overline{i})\right)\right)\nonumber\\
&&\phantom{t\left(\theta\right)=}\sum_{i}\left((ii)\otimes\left(d_{ii}(ii)+a_{ii}(\overline{i}\overline{i})\right)+
(\overline{i}\overline{i})\otimes\left(a_{ii}(ii)+d_{ii}(\overline{i}\overline{i})\right)\right)+\nonumber\\
&&\phantom{t\left(\theta\right)=}(pp)\otimes(pp)
\end{eqnarray}

Crucial features to be noted
\begin{enumerate}
    \item Only the diagonal blocks have nonzero terms on the diagonal;
    \item In each such blocks there are only two non-zero terms (with only one for the $p$-th);
    \item These features are iterated under successive coproducts.
\end{enumerate}
Thus
\begin{eqnarray}
&&\hbox{tr}\left(T\left(\theta\right)\right)=\hbox{tr}\left(\sum_i\left(t_{ii}\left(\theta\right)+
t_{\overline{i}\overline{i}}\left(\theta\right)\right)+t_{pp}\left(\theta\right)\right)\nonumber\\
&&\phantom{\hbox{tr}\left(T\left(\theta\right)\right)}=2\sum_i\left(d_{ii}+a_{ii}\right)+1\nonumber\\
&&\phantom{\hbox{tr}\left(T\left(\theta\right)\right)}=2\sum_ie^{m_{ii}^{(+)}\theta}+1\\
&&\hbox{tr}\left(T^{(2)}\left(\theta\right)\right)=2\sum_i\left(d_{ii}\left(d_{ii}+a_{ii}\right)+
a_{ii}\left(d_{ii}+a_{ii}\right)\right)+1\nonumber\\
&&\phantom{\hbox{tr}\left(T^{(2)}\left(\theta\right)\right)}=2\sum_i\left(d_{ii}+a_{ii}\right)^2+1\nonumber\\
&&\phantom{\hbox{tr}\left(T^{(2)}\left(\theta\right)\right)}=2\sum_ie^{2m_{ii}^{(+)}\theta}+1
\end{eqnarray}
and continuing stepwise
\begin{eqnarray}
&&\hbox{tr}\left(T^{(r)}\left(\theta\right)\right)=2\sum_{i=1}^{p-1}\left(d_{ii}+a_{ii}\right)^r+1
=2\sum_{i=1}^{p-1}e^{rm_{ii}^{(+)}\theta}+1
\end{eqnarray}
This is how the structure of our projector basis leads to a direct evaluation of $\hbox{tr}\left(T^{(r)}
\left(\theta\right)\right)$ for the general case without the full list of eigenvalues. The trace is given by
$2(p-1)+1=2p-1=N$ eigenvalues. The remaining $\left(N^r-N\right)$ eigenvalues give zero trace, as we have seen, due
to multiplet structures corresponding to roots of unity.

\section{\LARGE Encounter with a theorem of Fermat}
\setcounter{equation}{0}

A well-known theorem of Fermat states
\begin{equation}
N^r=N\,\hbox{mod.}\, r,
\end{equation}
where $\left(N,r\right)$ are positive integers and $r$ is a prime number. Writing it as
\begin{equation}
N^r-N=rM,
\end{equation}
we try to obtain the integer $M$ explicitly with the following purpose.

In sec. 6 we noted that out of $N^r$ eigenvalues of the transfer matrix $T^{(r)}\left(\theta\right)$ of order $r$
$\left(N^r-N\right)$ eigenvalues must give zero trace when summed. We have also seen how such a zero trace constraint is
implemented in our case through multiplets corresponding to roots of unity as explained in (6.8), (6.9) and (6.10).
When $r$ is a prime number the minimal multiplets can be consistently "r-plets" (or "$nr$-plets") only if (B.2) is satisfied. But precisely
this is guaranteed by (B.1). This is the link of our multiplet structure with (B.1). When $r$ is factorizable one can have
lower multiplets $\left(r_1,r_2,\ldots\right)$ for, say, $r=r_1\cdot r_2\cdots r_n$. This is already seen for $r=4$ in
(6.10). We are interested here in odd integers $N$, but (B.2) holds also for even $N$. We now construct $M$ giving the
number of $r$-plets.

Different constructions of $M$ are certainly possible. The one particularly suitable for our purpose is as follow. One has
for $r=1,\,3,\,5,\,7,\,11,\ldots$ respectively
\begin{eqnarray}
&&N-N=0,\nonumber\\
&&N^3-N=\left(N-1\right)N\left(N+1\right),\nonumber\\
&&N^5-N=\left(N-2\right)\left(N-1\right)N\left(N+1\right)\left(N+2\right)+5\left(N-1\right)N\left(N+1\right)\nonumber\\
&&N^7-N=\left(N-3\right)\left(N-2\right)\left(N-1\right)N\left(N+1\right)\left(N+2\right)\left(N+3\right)+\nonumber\\
&&\phantom{N^7-N=}7\left(2\left(N-2\right)\left(N-1\right)N\left(N+1\right)\left(N+2\right)+
3\left(N-1\right)N\left(N+1\right)\right)\nonumber
\end{eqnarray}
Continuing thus with product of consecutive factors
\begin{eqnarray}
&&N^{11}-N=\left(N-5\right)\cdots N\cdots\left(N+5\right)+\nonumber\\
&&\phantom{N^7-N=}11\left[5\left(N-4\right)\cdots N\cdots\left(N+4\right)+57\left(N-3\right)\cdots N\cdots\left(N+3\right)
+\right.\nonumber\\
&&\phantom{N^7-N=}128\left(N-2\right)\cdots N\cdots\left(N+2\right)+31\left(N-1N\left(N+1\right)\right]
\end{eqnarray}
Thus for $N=3$, $r=3,\,5,\,7,\ldots$, one has respectively 8 triplets, 48 5-plets, 312 7-plets and so on
(unless 10-plets, 14-plets and so on are also obtained).

The first term is
\begin{equation}
\left(N-\frac{r-1}2\right)\ldots \left(N+\frac{r-1}2\right)
\end{equation}
and being a product of $r$ consecutive integers, evidently divisible by $r$. The lower order products all have $r$ as
a factor. Hence the result. When $N\leq \frac{r-1}2$ one or more higher order products vanish. But through now $N^r$ is
not directly, visibly present on the right, the results still hold. Thus though ($N-3$) vanishes in the leading term,
\begin{equation}
3^7-3=0+7\cdot (312).
\end{equation}
For completeness we give the general result below
\begin{eqnarray}
&&N^r-N=\sum_{p=1}^{\frac{r-1}2}A_p\left(r\right)\prod_{k=-p}^p\left(N+k\right)\equiv
\sum_{p=0}^{\frac{r-1}2}A_p\left(r\right)B\left(N,p\right)
\end{eqnarray}
where
\begin{eqnarray}
&&A_{\frac{r-1}2}\left(r\right)=1,\nonumber\\
&&A_{\frac{r-2k+1}2}\left(r\right)=\sum_{m_1=1}^{k-1}\sum_{m_2=1}^{m_1-1}\cdots
\sum_{m_{k-1}=1}^{m_{k-2}-1}\left[-\frac{H_{m_1}^{2k-2m_1}\left(r\right)}{\left(2k-2m_1\right)!}\right]\times
\left[-\frac{H_{m_2}^{2m_1-2m_2}\left(r\right)}{\left(2m_1-2m_2\right)!}\right]\times\cdots\times\nonumber\\
&&\phantom{A_{\frac{r-2k+1}2}\left(r\right)=}\left[-\frac{H_{m_{k-1}}^{2m_{k-2}-2m_{k-1}}\left(r\right)}
{\left(2m_{k-2}-2m_{k-1}\right)!}\right],\qquad
k=2,3,\ldots,\frac{r-1}2.
\end{eqnarray}
The elements $H_k^m$ are given by
\begin{eqnarray}
&&H_k^m\left(r\right)=\left(\sum_{p_1=-\frac{r-2k+1}2}^{\frac{r-2k+1}2}p_1\right)
\left(\sum_{p_2=-\frac{r-2k+1}2\atop p_2\neq p_1}^{\frac{r-2k+1}2}p_2\right)\cdots
\left(\sum_{{p_m=-\frac{r-2k+1}2\atop p_m\neq p_1,\ldots,p_{m-1}}}^{\frac{r-2k+1}2}p_m\right)\qquad m\neq
0\nonumber\\
&&H_k^0\left(r\right)=1.
\end{eqnarray}
For example,
\begin{eqnarray}
&&H_k^2\left(r\right)=-\frac{1}{12}\left(r-2k+1\right)\left(r-2k+2\right)\left(r-2k+3\right),\nonumber\\
&&H_k^4\left(r\right)=\frac{1}{240}\left(5r+17-10k
\right)\left(r+3-2k\right)\left(r-2k+2\right)\left(r-2k+1\right)\left(
r-2k\right)\left(r-1-2k\right),\nonumber\\
&&H_{\frac{r-2k+1}2}^{2k}\left(r\right)=\left(-1\right)^k(2k)!(k!)^2,\,\,
k=0,\ldots,\frac{r-1}2,\nonumber\\
&&H_{\frac{r-2p+1}2}^{2p+2k-r-1}\left(r\right)=\left(2p+2k-r-1\right)!\left.\frac{d^{r-2k+2}B(N,p)}{dN^{r-2k+2}}
\right|_{N=0},\,\, p=\frac{{r-2k+1}}2,\ldots,\frac{{r-1}}2,\nonumber\\
&&H_{m_1}^{2k-2m_1}\left(r\right)=\left(2k-2m_1\right)!\left.\frac{d^{2m_1+1}B(N,\frac{r-2m_1+1}2)}{dN^{2m_1+1}}
\right|_{N=0},\,\, m_1=k,\ldots,\frac{{r-1}}2
\end{eqnarray}
The above formula gives the general expression for the coefficients. But, as Alain Lascoux has pointed out, the symmetric
form of Newton's interpolation fromula relevant for our case leads to complete functions as coefficients. These can be
obtained systematically and conveniently. Thus one obtains, for example,
\begin{eqnarray}
&&N^{11}-N=\left(N-5\right)\cdots N\cdots\left(N+5\right)+\nonumber\\
&&\phantom{N^7-N=}\left(N-4\right)\cdots N\cdots\left(N+4\right)\left(1^2+2^2+3^2+4^2+5^2\right)+\nonumber\\
&&\phantom{N^7-N=}\left(N-3\right)\cdots N\cdots\left(N+3\right)\left(1^4+2^4+3^4+4^4+1^2\cdot 2^2+1^2\cdot 3^2\right.\nonumber\\
&&\phantom{N^7-N=}\left.+1^2\cdot 4^2+2^2\cdot 3^2+2^2\cdot 4^2+3^2\cdot 4^2\right)+\nonumber\\
&&\phantom{N^7-N=}\left(N-2\right)\cdots N\cdots\left(N+2\right)\left(1^6+2^6+3^6+1^4\cdot 3^2+2^4\cdot 3^2\right.
\nonumber\\
&&\phantom{N^7-N=}\left.
+2^4\cdot 1^2+3^4\cdot 1^2+3^4\cdot 2^2+1^4\cdot 2^2+1^2\cdot 2^2\cdot 3^2\right)+\nonumber\\
&&\phantom{N^7-N=}\left(N-1\right)N\left(N+1\right)\left(1^8+2^8+1^6\cdot 2^2+2^6\cdot 1^2
+1^4\cdot 2^4\right)
\end{eqnarray}
This gives (B.3) in which 11 is already factorized.

\section{\LARGE $\widehat{R}tt$-algebra}
\setcounter{equation}{0}

We present below a canonical formulation of the $\widehat{R}tt$ algebra \cite{14} specifically adapted to our case. The
Baxterized form with $N^2$ blocks of $N\times N$ matrices $t_{ij}\left(\theta\right)$  must satisfy
\begin{equation}
\widehat{R}\left(\theta-\theta'\right)\left(t\left(\theta\right)\otimes t\left(\theta'\right)\right)=
\left(t\left(\theta'\right)\otimes t\left(\theta\right)\right)\widehat{R}\left(\theta-\theta'\right),
\end{equation}
where (since $P^2=I$)
\begin{equation}
t\left(\theta\right)\otimes t\left(\theta'\right)=
\left(t\left(\theta\right)\otimes I\right)\left(I\otimes t\left(\theta\right)\right)\equiv
t_1\left(\theta\right)t_2\left(\theta\right)=
\left(t_1\left(\theta\right)P\right)\left(Pt_2\left(\theta\right)\right).
\end{equation}
But $t_1\left(\theta\right)=Pt_2\left(\theta\right)P$, hence
\begin{equation}
t_1\left(\theta\right)P=t_2\left(\theta\right)P\equiv \widehat{t}\left(\theta\right).
\end{equation}
Thus
\begin{equation}
\widehat{R}\left(\theta-\theta'\right)\widehat{t}\left(\theta\right)\widehat{t}\left(\theta'\right)=
\widehat{t}\left(\theta'\right)\widehat{t}\left(\theta\right)\widehat{R}\left(\theta-\theta'\right),
\end{equation}
where one has just matrix multiplication of the same matrix $\widehat{t}$ with arguments $\left(\theta,\theta'\right)$.

Now suppose that one has obtained explicitly the diagonalizer $M$ of $\widehat{R}\left(\theta\right)$. When
\begin{equation}
\widehat{R}\left(\theta\right)=\sum_{\alpha,\beta}f_{\alpha\beta}\left(\theta\right)P_{\alpha\beta},
\end{equation}
where $P_{\alpha\beta}$ form a complete basis of $\theta$-independent projectors (and the minimal polynomial
equation satisfied by $\widehat{R}\left(\theta\right)$ has no multiple roots for consistency) one can construct a
$\theta$-independent $M$ to diagonalize each $P_{\alpha\beta}$ simultaneously. For our nested sequence of projectors
(1.3) the diagonalizer is given by in sec. 3 of ref. 1 as
\begin{eqnarray}
&&\sqrt{2}M=\sqrt{2}M^{-1}=\sqrt{2}(pp)\otimes (pp)+\nonumber\\
&&\phantom{\sqrt{2}M=\sqrt{2}M^{-1}=}(pp)\otimes\left(\sum_i\left( (ii)-(\overline{i}\overline{i})+(i\overline{i})+
(\overline{i}i)\right)\right)\nonumber\\
&&\phantom{\sqrt{2}M=\sqrt{2}M^{-1}=}\left(\sum_i\left( (ii)-(\overline{i}\overline{i})+(i\overline{i})+
(\overline{i}i)\right)\right)\otimes(pp)+\\
&&\phantom{\sqrt{2}M=\sqrt{2}M^{-1}=}\sum_{i,j}\left(\left((ii)-(\overline{i}\overline{i})\right)\otimes
\left((jj)+(\overline{j}\overline{j})\right)+\left((i\overline{i})+(\overline{i}i)\right)\otimes
\left((j\overline{j})+(\overline{j}j)\right)\right),\nonumber
\end{eqnarray}
where
\begin{equation}
i=1,2,\ldots,p-1,\qquad \overline{i}=2p-1,\ldots,p+1,\qquad N=2p-1.
\end{equation}
For $N=3$ ($p=2$), one obtains (see (3.5) of ref. 1)
\begin{eqnarray}
&&\sqrt{2}M=\sqrt{2}M^{-1}=\left(1,1,1,1,\sqrt{2},-1,-1,-1,-1\right)_{\hbox{diag.}}+\nonumber\\
&&\phantom{\sqrt{2}M=\sqrt{2}M^{-1}=}\left(1,1,1,1,\sqrt{2},1,1,1,1\right)_{\hbox{anti-diag.}}
\end{eqnarray}
($\sqrt{2}$ being the common element of diag and anti-diag).

The general case is now evident. Defining
\begin{equation}
M\widehat{R}\left(\theta\right)M^{-1}=D\left(\theta\right)
\end{equation}
a diagonal matrix and
\begin{equation}
M\widehat{t}\left(\theta\right)M^{-1}=K\left(\theta\right)
\end{equation}
quite generally
\begin{equation}
D\left(\theta-\theta'\right)K\left(\theta\right)K\left(\theta'\right)=
K\left(\theta'\right)K\left(\theta\right)D\left(\theta-\theta'\right).
\end{equation}
This our canonical formulation \cite{15}. For our present case (3.3), (3.4) of ref. \cite{1})
\begin{equation}
D\left(\theta\right)=\left(e^{m_{11}^{(+)}\theta},e^{m_{12}^{(+)}\theta},\ldots,e^{m_{11}^{(-)}\theta}\right).
\end{equation}
From (C.11)
\begin{equation}
D_{aa}\left(\theta-\theta'\right)K_{ac}\left(\theta\right)K_{cb}\left(\theta'\right)=
K_{ac}\left(\theta'\right)K_{cb}\left(\theta\right)D_{bb}\left(\theta-\theta'\right).
\end{equation}
For $N=3$, for our case, defining the $3\times 3$ diagonal blocks
\begin{eqnarray}
&&d_{11}\left(\theta\right)=\left(e^{m_{11}^{(+)}\theta},e^{m_{12}^{(+)}\theta},
e^{m_{11}^{(+)}\theta}\right)_{\hbox{diag;}}\nonumber\\
&&d_{22}\left(\theta\right)=\left(e^{m_{21}^{(+)}\theta},1,e^{m_{21}^{(-)}\theta}\right)_{\hbox{diag;}}\nonumber\\
&&d_{\overline{1}\overline{1}}\left(\theta\right)=\left(e^{m_{11}^{(-)}\theta},e^{m_{12}^{(-)}\theta},
e^{m_{11}^{(-)}\theta}\right)_{\hbox{diag;}}
\end{eqnarray}
and denoting $\left(K\left(\theta\right),K\left(\theta'\right),D\left(\theta-\theta'\right)\right)\equiv
\left(K,K',D''\right)$ when $d_{11}\left(\theta-\theta'\right)\equiv d_{11}''$ and so on
\begin{equation}
\left(\begin{array}{ccc}
  d_{11}''\left(KK'\right)_{11} & d_{11}''\left(KK'\right)_{12} & d_{11}''\left(KK'\right)_{1\overline{1}} \\
  d_{22}''\left(KK'\right)_{21} & d_{22}''\left(KK'\right)_{22} & d_{22}''\left(KK'\right)_{2\overline{1}} \\
  d_{\overline{1}\overline{1}}''\left(KK'\right)_{\overline{1}1} & d_{\overline{1}\overline{1}}''
  \left(KK'\right)_{\overline{1}2} & d_{\overline{1}\overline{1}}''\left(KK'\right)_{\overline{1}\overline{1}} \\
\end{array}\right)=\left(\begin{array}{ccc}
  \left(K'K\right)_{11}d_{11}'' & \left(K'K\right)_{12}d_{22}'' & \left(K'K\right)_{1\overline{1}}
  d_{\overline{1}\overline{1}}'' \\
  \left(K'K\right)_{21}d_{11}'' & \left(K'K\right)_{22}d_{22}'' & \left(K'K\right)_{2\overline{1}}
  d_{\overline{1}\overline{1}}'' \\
  \left(K'K\right)_{\overline{1}1}d_{11}'' & \left(K'K\right)_{\overline{1}2}d_{22}''
  & \left(K'K\right)_{\overline{1}\overline{1}}d_{\overline{1}\overline{1}}'' \\
\end{array}\right).
\end{equation}
Substituting the explicit form of $\left(KK'\right)_{ij}$ one obtains the full set of 81 relations (for $N=3$) of the
$\widehat{R}tt$-algebra. Each $K_{ij}$ is a $3\times 3$ block whose elements are the blocks of $t_{ij}$. Thus
\begin{eqnarray}
&&2K_{11}=\left(\begin{array}{ccc}
  t_{11}+t_{\overline{1}\overline{1}} & t_{12}+t_{\overline{1}2} & t_{1\overline{1}}+t_{\overline{1}1} \\
  0 & 0 & 0 \\
  t_{1\overline{1}}+t_{\overline{1}1} & t_{12}+t_{\overline{1}2} & t_{11}+t_{\overline{1}\overline{1}} \\
\end{array}\right)\nonumber\\
&&\phantom{2K_{11}}=\left(t_{11}+t_{\overline{1}\overline{1}}\right)\left((11)+(\overline{1}\overline{1})\right)
+\left(t_{12}+t_{\overline{1}2}\right)\left((12)+(\overline{1}2)\right)+\nonumber\\
&&\phantom{2K_{11}=}\left(t_{1\overline{1}}+t_{\overline{1}1}\right)\left((1\overline{1})+(\overline{1}1)\right).
\end{eqnarray}
In such a notation
\begin{eqnarray}
&&2K_{\overline{1}\overline{1}}=\left(t_{11}+t_{\overline{1}\overline{1}}\right)\left(-(11)+(\overline{1}\overline{1})
\right)+\left(t_{12}+t_{\overline{1}2}\right)\left(-(12)+(\overline{1}2)\right)+\nonumber\\
&&\phantom{2K_{11}=}\left(t_{1\overline{1}}+t_{\overline{1}1}\right)\left(-(1\overline{1})+(\overline{1}1)\right),\\
&&2K_{1\overline{1}}=\left(t_{11}-t_{\overline{1}1}\right)\left((11)-(\overline{1}\overline{1})
\right)+\left(t_{12}-t_{\overline{1}2}\right)\left((12)-(\overline{1}2)\right)+\nonumber\\
&&\phantom{2K_{11}=}\left(t_{11}-t_{\overline{1}\overline{1}}\right)\left((1\overline{1})-(\overline{1}1)\right),\\
&&2K_{\overline{1}1}=\left(t_{1\overline{1}}-t_{\overline{1}1}\right)\left((11)+(\overline{1}\overline{1})
\right)+\left(t_{12}-t_{\overline{1}2}\right)\left((12)+(\overline{1}2)\right)+\nonumber\\
&&\phantom{2K_{11}=}\left(t_{11}-t_{\overline{1}\overline{1}}\right)\left((1\overline{1})+(\overline{1}1)\right),\\
&&2K_{12}=\left(t_{11}+t_{1\overline{1}}+t_{\overline{1}1}+t_{\overline{1}\overline{1}}\right)(21)+
\sqrt{2}\left(t_{12}+t_{\overline{1}2}\right)(22)\nonumber\\
&&\phantom{2K_{11}=}\left(t_{11}-t_{\overline{1}\overline{1}}+t_{\overline{1}1}-t_{\overline{1}\overline{1}}
\right)(2\overline{1}),\\
&&2K_{\overline{1}2}=\left(t_{11}+t_{1\overline{1}}-t_{\overline{1}1}-t_{\overline{1}\overline{1}}\right)(21)+
\sqrt{2}\left(t_{12}-t_{\overline{1}2}\right)(22)\nonumber\\
&&\phantom{2K_{11}=}\left(t_{11}-t_{\overline{1}\overline{1}}-t_{\overline{1}1}+t_{\overline{1}\overline{1}}
\right)(2\overline{1}),\\
&&2K_{21}=\left(t_{21}+t_{2\overline{1}}\right)\left((11)+(\overline{1}\overline{1})\right)+2t_{22}(12)
+\left(t_{21}-t_{2\overline{1}}\right)\left((\overline{1}1)-(\overline{1}\overline{1})\right),\\
&&2K_{2\overline{1}}=\left(t_{21}-t_{2\overline{1}}\right)\left(-(11)+(1\overline{1})\right)+2t_{22}(\overline{1}2)
+\left(t_{21}+t_{2\overline{1}}\right)\left((\overline{1}1)+(\overline{1}\overline{1})\right),\\
&&2K_{22}=\sqrt{2}\left(t_{21}+t_{2\overline{1}}\right)(21)+
2t_{22}(22)+\sqrt{2}\left(t_{21}-t_{2\overline{1}}\right)(2\overline{1}).
\end{eqnarray}
We have not directly utilized the $\widehat{R}tt$ constraints in constructing eigenstates. But since (C.1) is the basic
equation providing the starting point we present here the most compact approach to the full set of 81 constraints for
$N=3$.

\end{appendix}

\end{document}